\documentstyle[11pt]{article}
\newtheorem{corr}{Corollary}[section]

\newtheorem{theo}{Theorem}[section]
\newtheorem{lemma}{Lemma}[section]
\newtheorem{defi}{Definition}[section]
\newtheorem{rem}{Remark}[section]

\newtheorem{con}{Conjecture}[section]
\newcommand{\im}{{\rm im}}
\newcounter{bean}
\newcounter{inner}
\begin{document}
\addtolength{\baselineskip}{4pt}
\makeatletter
%
\expandafter\ifx\csname amssym.def\endcsname\relax \else\endinput\fi
%
\expandafter\edef\csname amssym.def\endcsname{%
       \catcode`\noexpand\@=\the\catcode`\@\space}
\catcode`\@=11
%

\def\undefine#1{\let#1\undefined}
\def\newsymbol#1#2#3#4#5{\let\next@\relax
 \ifnum#2=\@ne\let\next@\msafam@\else
 \ifnum#2=\tw@\let\next@\msbfam@\fi\fi
 \mathchardef#1="#3\next@#4#5}
\def\mathhexbox@#1#2#3{\relax
 \ifmmode\mathpalette{}{\m@th\mathchar"#1#2#3}%
 \else\leavevmode\hbox{$\m@th\mathchar"#1#2#3$}\fi}
\def\hexnumber@#1{\ifcase#1 0\or 1\or 2\or 3\or 4\or 5\or 6\or 7\or
8\or
 9\or A\or B\or C\or D\or E\or F\fi}

\font\tenmsa=msam10
\font\sevenmsa=msam7
\font\fivemsa=msam5
\newfam\msafam
\textfont\msafam=\tenmsa
\scriptfont\msafam=\sevenmsa
\scriptscriptfont\msafam=\fivemsa
\edef\msafam@{\hexnumber@\msafam}
\mathchardef\dabar@"0\msafam@39
\def\dashrightarrow{\mathrel{\dabar@\dabar@\mathchar"0\msafam@4B}}
\def\dashleftarrow{\mathrel{\mathchar"0\msafam@4C\dabar@\dabar@}}
\let\dasharrow\dashrightarrow
\def\ulcorner{\delimiter"4\msafam@70\msafam@70 }
\def\urcorner{\delimiter"5\msafam@71\msafam@71 }
\def\llcorner{\delimiter"4\msafam@78\msafam@78 }
\def\lrcorner{\delimiter"5\msafam@79\msafam@79 }
\def\yen{{\mathhexbox@\msafam@55 }}
\def\checkmark{{\mathhexbox@\msafam@58 }}
\def\circledR{{\mathhexbox@\msafam@72 }}
\def\maltese{{\mathhexbox@\msafam@7A }}

\font\tenmsb=msbm10
\font\sevenmsb=msbm7
\font\fivemsb=msbm5
\newfam\msbfam
\textfont\msbfam=\tenmsb
\scriptfont\msbfam=\sevenmsb
\scriptscriptfont\msbfam=\fivemsb
\edef\msbfam@{\hexnumber@\msbfam}
\def\Bbb#1{{\fam\msbfam\relax#1}}
\def\widehat#1{\setbox\z@\hbox{$\m@th#1$}%
 \ifdim\wd\z@>\tw@ em\mathaccent"0\msbfam@5B{#1}%
 \else\mathaccent"0362{#1}\fi}
\def\widetilde#1{\setbox\z@\hbox{$\m@th#1$}%
 \ifdim\wd\z@>\tw@ em\mathaccent"0\msbfam@5D{#1}%
 \else\mathaccent"0365{#1}\fi}
\font\teneufm=eufm10
\font\seveneufm=eufm7
\font\fiveeufm=eufm5
\newfam\eufmfam
\textfont\eufmfam=\teneufm
\scriptfont\eufmfam=\seveneufm
\scriptscriptfont\eufmfam=\fiveeufm
\def\frak#1{{\fam\eufmfam\relax#1}}
\let\goth\frak

\csname amssym.def\endcsname
\makeatother

\title{Factorization  theorems for the  representations of the fundamental
groups of quasiprojective varieties and some applications }
\author{L. Katzarkov}
\maketitle

\begin{abstract} In this paper,  using Gromov-Jost-Korevaar-Schoen technique
of harmonic maps to nonpositively curved targets,  we study  the
representations of the fundamental groups of quasiprojective varieties. As an
application of the above considerations we give a proof of a weak version of
the Shafarevich Conjecture.
\end{abstract}

\section {Introduction}

It is known  \cite{10} that to every representation of the fundamental group of
 a smooth projective variety $X$ to a simple complex Lie  group $G$  we can
assign  an equivariant  harmonic map from the universal cover of $X$ to the
corresponding symmetric space. If $G$ is a simple Lie group over an arbitrary
local field $K$, then we can try to work with equivariant  harmonic maps from
the universal cover of $X$ to the corresponding Euclidean building $B$ instead.

Developing the theory of the equivariant  harmonic maps from the universal
cover of a  Riemannian manifold $X$  to  $B$, a Euclidean building, Gromov and
Schoen  in \cite{1} and Korevaar and Schoen in  \cite{2}  proved the following
result:

\begin{theo} Consider a Zariski  dense representation of the fundamental group
of $X$, $\pi_{1}(X)$, to some simple Lie group $G$  over $K$. Then if
$\pi_{1}(X)$  acts on the Euclidean building  $B$ corresponding to $G$ without
fixed points, there exists an equivariant  nonconstant harmonic map $U$ from
the universal cover of the  Riemannian manifold $X$  to  $B$.
\end{theo}

The main goal  of this paper is to study the factorizations of Zariski  dense
representations of the algebraic fundamental group of $X$ to some simple Lie
group $G$  over $K$ using the theory of harmonic maps. The pioneer ideas in
this direction belong to Eells and Sampson, Carlson and Toledo and  Yau and
Jost  (see \cite{Y} for more detailed account on that). Later  these ideas were
developed by  Corlette and Simpson \cite{12}, \cite{11} and Mok \cite{20}, who
uses the existence of a K\"{a}hler form on a symmetric space. Since there is no
K\"{a}hler form on the building $B$ we cannot use Mok's technique \cite{20} to
get a holomorphic foliation on $X$. The Higgs bundles technique from \cite{2}
and \cite{3} also  does not work. Instead we use a generalization of the
Clemens-Lefschetz-Simpson theorem from \cite{3}  and the Gromov-Schoen theory
to prove the   main theorem:

\begin{theo} Let $\varrho:\pi_{1}(X) \longrightarrow G$ be a
 Zariski  dense representation of the   fundamental group of $X$ to some simple
Lie group $G$  over $K$.  Then
\begin{list}{{\bf \Alph{bean}.}}{\usecounter{bean}}
\item either  the image of $\varrho$ is  in a maximal compact subgroup of $G$
or,
\item there exist:
\begin{list}{\arabic{inner})}{\usecounter{inner}}
\item  a finite etale cover $X'$ of a blow up of $X$;
\item a smooth projective variety $Y$ of positive dimension    $l$ less then or
equal to the rank $r$  of $G$ over $K$;
\item  a holomorphic map $h : X' \longrightarrow Y $
such that $\varrho:\pi_{1}(X') \longrightarrow G$ factors through a
representation of $\pi_{1}(Y) $.
\end{list}
\end{list}
\end{theo}

We will always assume that the representation $\varrho :
\pi_{1}(X) \longrightarrow G$ is centerless, namely the intersection
of  $\varrho(\pi_{1}(X))$ and the center of $G$, $Z(G)$ is trivial.

We also partially  extend this result to the case when $X \setminus D$ is a
quasiprojective variety and $\varrho : \pi_{1}(X) \longrightarrow G$ is a
representation  with  unipotent monodromy at infinity.

Recently Simpson  and Corlette have proven the following theorem:

\begin{theo} Let  $S$ be a smooth projective variety and $\varrho$ be a Zariski
dense representation of $\pi_{1}(S)$ in $SL(2,{\Bbb C})$. Then one of the
following is true:
 \begin{list}{\arabic{inner})}{\usecounter{inner}}
 \item $\varrho$ is a pullback from an orbicurve;
 \item $\varrho$ is a pullback from a Hilbert modular variety.
\end{list}
\end{theo}

This result was partially generalized in  \cite{2} to higher rank complex
groups.

Let $\varrho:\pi_{1}(X) \longrightarrow G$ be a  Zariski  dense nonrigid
representation of the   fundamental group of $X$ to some  complex  simple Lie
group $G$.  Then there exists a  curve $I$ in the moduli space of
representations passing through
 $\varrho :\pi_{1}(X) \longrightarrow G$ in direction of which $\rho$ deforms
nontrivialy. Since the moduli space of representations is an affine variety the
curve $I$  is affine as well. Let $\bar{I}$ be the compactification of $I$.
Consider now the representation  $\bar{\rho}$ to $\bar{G}$, where $ \bar{G} $
is defined over the field of  fractions of  $I$. Let $O_{p}$ be the
localization at some point $p \in \bar{I} \setminus I$. Let $Z[T]$ be an
extension of the ring of integers of  $O_{p}$ which contains all coefficients
of $\bar{\rho}$ localized at the point $p$. Consider now $I$ as a curve over
${\rm Spec} Z[T]$. The fact that $\varrho
: \pi_{1}(X) \longrightarrow G$ is a   Zariski  dense nonrigid representation
implies that $\bar{\rho}$ is   a   Zariski  dense nonrigid representation.
Observe that there is still going to be a point $q$ in the  curve $\bar{I}$
which is not in ${\rm Spec}
Z[T]$. Let $\chi(\bar{\rho})$ be a character of  $\bar{\rho}$. Moving over
${\rm Spec} Z[T]$ in  $\bar{I}$ we see that $\chi(\bar{\rho})$ is unbounded at
$q$. Therefore we conclude that  the  representation $\bar{\rho}$ to $\bar{G}$
is not contained in any bounded subgroup in $\bar{G}$.

Using Theorem 1.2 we were able to prove the following:

\bigskip

\noindent
{\bf Corollary 4.2} {\it Let $\varrho:\pi_{1}(X) \longrightarrow G$ be a
Zariski  dense nonrigid representation of the   fundamental group of $X$ to
some  complex  simple Lie group $G$. Then  there exist:
\begin{list}{\arabic{inner})}{\usecounter{inner}}
\item  a finite etale cover $X'$ of a blow up of $X$;
\item a smooth projective variety $Y$ of positive dimension  $l$
less or equal to the rank $r$  of $G$ over $\Bbb{C}$;
\item  and a  holomorphic map $h : X' \longrightarrow Y$
such that $\varrho : \pi_{1}(X') \longrightarrow G$ factors through a
representation of $\pi_{1}(Y)$.
\end{list}}

 \bigskip

In the complimentary situation, namely when we consider a Zariski  dense rigid
representation of the   fundamental group of $X$ to a complex  simple Lie group
$G$, following Simpson \cite{3} we can assign to every Zariski dense rigid
representation  a new  Zariski dense rigid representation $\varrho: \pi_{1}(X)
\longrightarrow G^{1}$  into a group $G^{1}$ defined over a local field. The
procedure goes as follows.

Observe that the  moduli space of representations is defined over ${\Bbb Q}$,
and since we are working with a rigid representation
 we can find an isomorphic representation defined over $\bar{\Bbb
 Q}$. Therefore  we can assume that $K = \bar{\Bbb Q}$.
Let $E$ be a finite extension of  ${\Bbb Q}$ defined to  be the extension
which contains all coefficients of our representation, and let $O$ denote the
ring of integers in $E$.
Let $E_{p}$ denote the field of fractions of the completion of $O$  at $p$, for
some  prime $p$. Let $G^{1}$ be the  group
of  $E_{p}$ valued points of $G$ and we again  use $\varrho$ for the
representation $\varrho: \pi_{1}(X) \longrightarrow G^{1}$. Therefore we get:

\bigskip

\noindent
{\bf Corollary 7.1} {\it Let $\varrho: \pi_{1}(X) \longrightarrow G$ be a
Zariski dense rigid representation. Then one of the following holds :
\begin{list}{{\bf \Alph{bean}.}}{\usecounter{bean}}
\item For every prime $p$ the image of $\varrho: \pi_{1}(X) \longrightarrow
G^{1}$ is contained in a maximal compact subgroup in $G^{1}$.
\item For   some prime $p$ the image of $\varrho: \pi_{1}(X) \longrightarrow
G^{1}$ is not  contained in a  maximal compact subgroup in $G^{1}$. Then there
exist:
\begin{list}{\arabic{inner})}{\usecounter{inner}}
\item  a finite etale cover $X'$ of a blow up of $X$;
\item  a smooth projective variety $Y$ of positive dimension  $l$ less or equal
to the rank $r$  of $G$ over $\Bbb{C}$;
\item and a  holomorphic map $h : X' \longrightarrow Y$
such that $\varrho : \pi_{1}(X') \longrightarrow G$ factors through a
representation of $\pi_{1}(Y)$.
\end{list}
\end{list} }

\bigskip

The paper is organized as follows:

The first part of the paper contains 4 sections.
After the introduction in section 2 we prove the
Clemens-Lefschetz-Simpson theorem which generalizes a previous result by
Simpson and is  a main tool in our proofs.  We give a brief explanation  of the
Gromov-Korevaar-Schoen-Jost  theory in  section 3. The proof of  our main
theorem is given in section 4.

We consider  some applications of the above  results in the second  part of the
paper. Section 5 contains the proof of the main  theorem in the quasiprojective
case.

Following the  general philosophy of  studying  the fundamental groups  of
projective varieties  by studying their representations,   we consider   in
section 6 the representations of the  fundamental  groups of projective
varieties onto the fundamental  groups of  some negatively curved polyhedra
defined by Benakli in \cite{32}. Using  the Clemens-Lefschetz-Simpson theorem
we give a partial  answer to a  question posed by Gromov in \cite{22}, proving
the following theorem for the Benakli's complexes:

\bigskip

\noindent
{\bf Theorem 6.1} {\it Every representation of the  fundamental  group of a
projective variety onto the fundamental  group of  some negatively curved
2-dimensional polyhedra, comes from the representation of the fundamental group
of an orbicurve.}

\bigskip

In section 7 using a recent work of N. Katz  \cite{36} we prove some
generalizations  of the  result of Corlette and Simpson. We say that a
representation is  integral if it is conjugate to a representation whose matrix
coefficients  are algebraic integers and that a representation is of dimension
one  if the corresponding harmonic map to the building $B$ is of rank one.

We were not able to get further then   the following:

\bigskip

\noindent
{\bf Theorem 7.1} {\it Let $\varrho: \pi_{1}(X) \longrightarrow G$ be  a rigid
Zariski dense representation of dimension  one. Then it is integral, in other
words it is conjugate to a representation whose matrix coefficients  are
algebraic integers.}
\bigskip

The idea is that since the algebraic curves do not have many rigid
representations something which factors through them does not have either.

In the last  section we make a connection with the theory of  Shafarevich maps

\[ Sh : X \dashrightarrow Sh(X)\]
developed by J\'anos Koll\'ar \cite{30}, \cite{31}. We consider also the
relative version of this map

\[Sh^{H} : X \dashrightarrow Sh^{H}(X),\]
 where $H$ is a normal subgroup in $\pi_{1}(X)$.

Using Theorems 1.2 and 1.3 we were able to estimate the dimension of the image
of the  Shafarevich map:

\bigskip

\noindent
{\bf Theorem 8.2} {\it Let X be a smooth projective variety  which has a  type
{\bf B} (from Theorem 1.2) Zariski dense representation $\rho$ of its
fundamental group  to a Lie  group $G$ defined over a  local field $K$. Define
$H= \ker (\rho)$. Then:

\[ \dim_{\Bbb{C}}Sh^{H}(X) = \dim_{\Bbb{C}}(Y)\leq {\rm
rank}_{\Bbb{C}}(G),\]

where $Y$ is the variety defined in Theorem 1.2. Moreover there exists a
finite nonramified covering of $Y$ which is  birationally isomorphic to a
finite nonramified covering of $Sh^{H}(X)$.}

Using the above theorem we are able to show that  in some cases the Shafarevich
morphism (in the sense of J. Koll\'ar)  exists. Namely, we prove the following
theorem:

\bigskip

\noindent
{\bf Theorem 8.3} {\it Let X be a smooth projective variety  which has a  type
{\bf B} (from Theorem 1.2) Zariski dense representation $\rho$ of its
fundamental group  to a Lie  group $G$ defined over a  local field $K$. Define
$H=\ker (\rho)$. Then: $Y$ is isomorphic to $Sh^{H}(X')$ and the map

\[{\bf Sh}^{H} : X' \longrightarrow {\bf Sh}^{H}(X'),\]
is a morphism. Here $ X'$ and $Y$ are the same as in Theorem 1.2.}

\bigskip

The Shafarevich conjecture says that for every smooth projective variety $X$
there exists a  Stein manifold ${\bf Sh}(\widetilde{X} )$ and a proper map with
connected fibers $ {\bf Sh} : \widetilde{X} \longrightarrow {\bf
Sh}(\widetilde{X})$. It was shown in  \cite{30} and  \cite{31} that the
following are equivalent:

1) There exists a  space  ${\bf Sh}(\widetilde{X} )$ with no compact complex
subspaces and a proper map with connected fibers $ {\bf Sh} : \widetilde{X}
\longrightarrow {\bf Sh}(\widetilde{X})$. This differs from the Shafarevich
conjecture since we have weaken the requirement about ${\bf Sh}(\widetilde{X}
)$ being Stein.

 2) There exists  a holomorphic map ${\bf Sh} : X \longrightarrow {\bf Sh}(X)
$, which contracts all subvarieties  $Z$
 in $X$   having   the property that $\im [\pi_{1}(Z)\longrightarrow
\pi_{1}(X)]$ is finite.

( We use  ${\bf Sh} : X \longrightarrow  {\bf Sh}(X) $ for the Shafarevich
morphism, instead of  $Sh : X \dashrightarrow Sh(X)$, which we use for the
rational map defined by  J\'anos Koll\'ar.) The existence of the moprphism

\[{\bf Sh} : X \longrightarrow {\bf Sh}(X)\]
for $X$ smooth projective was conjectured by J\'anos Koll\'ar in \cite{30} and
\cite{31}. This conjecture we call a weak version of the Shafarevich conjecture
or Shafarevich-Koll\'ar conjecture.

Therefore using theorem 8.3 we obtain the following:

\bigskip

\noindent
{\bf Corollary 8.1} {\it  If in the condition of the previous theorem we have
that $H = \ker(\rho)$ is finite then Shafarevich-Koll\'ar conjecture is true .}

\bigskip

The  Shafarevich-Koll\'ar conjecture in general is connected with another
question of Gromov:

\bigskip
\noindent
{\bf Question.} {\it
 Can we find a faithful  discrete cocompact action of every word hyperbolic
K\"{a}hler group on a space with $K<0$?}

\bigskip

The theorem of Rips  \cite{22} gives an answer to this question in some cases.
In
these cases, using the technique developed in this paper one can try to show
that  the Shafarevich-Koll\'{a}r conjecture follows for every  K\"{a}hler group
for which the map:

\[U : \widetilde{X} \longrightarrow T \]
is pluriharmonic.

We were able to realize the above idea  for fundamental groups, which have
faithful discrete cocompact action on trees and
2-dimensional negatively curved 2-complexes discussed in
section 6.

As a consequence we obtain the following nonvanishing theorem:

\bigskip

\noindent
{\bf Theorem 8.4}
{\it Let $X $ be a smooth fourfold of general type which has a  type {\bf B.}
( Theorem 1.2)  Zariski dense representation of its algebraic fundamental
group, $\rho: \pi_{1}(X) \longrightarrow
G$, where $G$ is a  simple Lie group of rank two or three over a local field
$K$. Then either:
\begin{list}{\arabic{inner})}{\usecounter{inner}}
\item   $\rho: \pi_{1}(X) \longrightarrow G $ factors through a representation
of the fundamental group of an orbicurve, or
\item $P_{n}(X) := H^{0} (X,  nK_{X})$ is not zero for $n \geq 4$.
\end{list}}
\bigskip

\noindent
{\bf Acknowledgements:}  I am  deeply indebted to  K. Corlette, P. Deligne, M.
Gromov, N. Katz, J. Koll\'ar, R. Schoen and C. Simpson for introducing  me  to
the subject and for their  constant interest in  my   work. I have also
benefited a lot  from  conversations with J. Block, N. Benakli,  J.L.
Brylinski, J. Carlson, M. Davis, R. Donagi, C. L. Epstein, W. Goldman, J. Jost,
M. Larsen, D.  McLaughlin, T. Pantev, M. Ramachandran, J. Shaneson, D. Toledo,
S. T. Yau and K. Zuo.

\begin{rem} {\rm Recently we were informed by J. Jost and K. Zuo that they have
obtained results similar to our results from  section 4.}
\end{rem}

\section{Lefschetz theorem}

In this section we generalize a result of Simpson \cite{3}. The
Clemens-Lefschetz-Simpson theorem we  prove  is of independent
interest but it is also one of the key arguments in  the next sections.

Let $X$ be a smooth projective variety and let $\alpha _{1},
\alpha _{2}, \ldots ,\alpha _{n}$ be linearly independent holomorphic one-forms
in $H^{0}(X, \Omega^{1}_{X})$. Denote by $\widetilde{X}$ the  universal
covering of $X$. Then by integrating  the pullbacks  of the forms  $\alpha
_{1}, \alpha _{2}, \ldots ,\alpha _{n}$ over $\widetilde{X}$ we can define a
map

\[ g: \widetilde{X} \longrightarrow {\Bbb C}^{n}.\]
 Let us first consider the case $n=2$.

\begin{theo} Let $\alpha _{1}, \alpha _{2}$  be two linearly
independent  holomorphic one forms.  Then one of the following cases holds:
\begin{list}{\alph{bean})}{\usecounter{bean}}
\item  The map $g$ has connected fibers.
\item There is a holomorphic map with connected fibers from $X$ to a normal
projective variety $Y$ such that $1 \leq \dim Y \leq 2$ and the forms $\alpha
_{1}, \alpha _{2}$  are pullbacks from $Y$.
\item The third case is a combination of the previous two. Namely, up to a
linear change of the coordinates in ${\Bbb C}^{2}$ the first form  $\alpha
_{1}$ is a pullback by a map   $X \longrightarrow Y$, where Y is an algebraic
curve.  The second form  $\alpha _{2}$ gives us a map  $ \widetilde{X}
\longrightarrow  {\Bbb C}$ with connected fibers.
\end{list}
\end{theo}
{\bf Proof.} We  borrow  heavily  from Simpson's arguments from  \cite{3},
where he proves the case $n=1$.

Let $Alb (X)$ be  the Albanese variety of $X$. The forms
$\alpha_{1}, \alpha_{2}$ are pullbacks of forms on $Alb (X)$, which we will
also call  $\alpha _{1}, \alpha_{2}$.
Define $B$ as the biggest abelian subvariety over which both forms  $\alpha
_{1}$ and  $\alpha_{2}$ considered as a sections of
$\Omega^{1}_{Alb (X)}$ are zero. Let

\[ A=Alb(X)/B.\]
We have the natural map

\[a : X \longrightarrow  A.\]
Let $l : \widetilde{X} \longrightarrow X$ be the universal covering map.

Suppose first that the dimension of $A$ is zero. But this means that the forms
$\alpha _{1}, \alpha_{2}$ are the zero forms.  Therefore $\dim_{\Bbb{C}}A \geq
1$.

Consider the map $g$ . Let $Q$ be the subset of ${\Bbb C^{2}}$ over which  the
map $g$ has multiple fibers. From the construction of $g$ it follows that $Q$
is a constructible set in ${\Bbb C^{2}}$.  We show now that if  $\dim_{\Bbb C}Q
> 0$ then  either case b) or c)  of the theorem holds.

Assume that $\dim_{\Bbb C}Q > 0$. Then we have that $Q$ is an analytic
subvariety of  ${\Bbb C}^{2}$ and  $\dim_{\Bbb C}Q = 1$.

We consider now the composition of the maps

\[g: \widetilde{X} \longrightarrow {\Bbb C}^{2} \longrightarrow {\Bbb C}^{1},\]
where the map

\[p : {\Bbb C}^{2} \longrightarrow {\Bbb C}^{1} = {\Bbb C}^{2} / l \]

is just a quotient by  some line  $l$ in ${\Bbb C^{2}}$.

Let us choose coordinates $t_{1}$ and $t_{2}$  in ${\Bbb C}^{2}$ so  that
$\alpha_{1}=g^{*}(dt_{1})$ and $\alpha_{2}=g^{*}(dt_{2})$.
There are two main cases to be dealt with.

Case $\alpha$) There is a  line

 \[l: c_{1}t_{1} +c_{2}t_{2}=0 ,\]

contained  in $Q$. Here $c_{1}$ and $c_{2}$ are some complex numbers.

In this case  the map

\[g: \widetilde{X} \longrightarrow {\Bbb C}^{2}  \longrightarrow  {\Bbb C}^{1}
\]

is given by integration of a holomorphic one form. Applying   Simpson's version
from \cite{3}  we obtain a map $X \longrightarrow Y$, where $Y$ is an algebraic
curve. This places us  in cases b) or c)  of the theorem. We will  be in case
b)  if there exists another line   in $Q$, which is transversal to   $l$.
Otherwise we get case c) of the theorem.

\medskip

Case $\beta$)  $Q$ does not contain a line.

Let $x$  be a smooth  point in $Q$. Then we need to consider two possibilities.

\begin{enumerate}
\item The differential $g_{*y}$ is the zero map on $T\widetilde{X}_{y}$ for
every $y$ from $g^{-1}(x)$.
We show that in this situation we have a map $X \longrightarrow Y$, where $Y$
is an algebraic surface. Consider the map
\[a : X \longrightarrow  A.\]
Let the map $f : X \longrightarrow Y$ be the Stein factorization of $a$.
Pulling back the map $f : X \longrightarrow Y$ to the
universal coverings of $X$ and $Y$, $\widetilde{X}$  and $\widetilde{Y}$, we
get that the map
\[ g : \widetilde{X} \longrightarrow  {\Bbb C}^{2}  \]
factors through the map
\[\tilde{f} : \widetilde{X} \longrightarrow \widetilde{Y}.\]
Suppose now that $\dim_{\Bbb{C}}Y>2$. Obviously a connected component of the
fiber of $\tilde{f}$ maps to a component of $v$, where $v$ is the map
\[v:\widetilde{Y} \longrightarrow {\Bbb C^{2}}, \]
with the property that $g=\tilde{f}\circ v$. Note that under the assumption
$\dim_{\Bbb{C}}Y > 2$ we have that dimension of the fiber of $v$ is at least
one.
But  the fact that $T\widetilde{X}_{y}$ for every $y$ from $g^{-1}(x)$ goes to
a point under $g_{*}$ means  that the restrictions and projections of the forms
$\alpha_{1}$ and $  \alpha_{2}$  on  a connected component of $g^{-1}(x)$ are
equal to zero. This implies that $g^{-1}(x)$ goes to a discrete set of points
under $ \widetilde{f}$, which contradicts  the fact that the dimension of the
fiber of $v$ is at least one. Therefore there exists a map $X \longrightarrow
Y$, where $Y$ is an algebraic surface, so we are in case b) of the theorem.
\item  For every $y$ from $g^{-1}(x)$, $g_{*}$ sends $T\widetilde{X}_{y}$  to
the same  line in
$T \times {\Bbb C^{2}}$. Let $N$ be the irreducible component of $Q$ passing
through $x$.
Since $Q$ does not contain a line  we can assume that the  two form
\[ dt_{1}  \wedge  dt_{2}\]
is zero when restricted to $N$ and  $ dt_{1},  dt_{2}$ are linearly independent
at every point $x$ in $N$ .
In this case  the Castelnuovo- De Francis theorem gives us  a  holomorphic map
$s : l(g^{-1}(N)) \longrightarrow Y$, where $Y$ is an algebraic curve. We  can
see that the map $g : g^{-1}(N) \longrightarrow \widetilde{Y}$ factors through
the map $\widetilde{s} : g^{-1}(N) \longrightarrow \widetilde{Y}$.
Consider now the map
\[a:X \longrightarrow  A.\]
Let $J=a(X)$. We prove next that $\dim_{\Bbb C} J \leq2$. Suppose not. From the
construction of $a$ we see that the map
\[g: \widetilde{X} \longrightarrow {\Bbb C}^{2} \]
factors through the maps
\[\widetilde{a}: \widetilde{X} \longrightarrow \widetilde{J},\]
\[v: \widetilde{J} \longrightarrow  {\Bbb C}^{2}.\]
To get a contradiction we need to show that there exists a fiber $V$ of the map
$v$ such that $\dim_{\Bbb C}V=0$.
We know that $\widetilde{Y} $ is contained in $\widetilde{J}$. But
by the universality of the Albanese  map for $l(g^{-1}(N))$ we can see that
the map
\[s: l(g^{-1}(N)) \longrightarrow Y \]
is nothing else but the map
\[a:X \longrightarrow  A,\]
restricted on $ l(g^{-1}(N)$. Recall that the map
\[g: \widetilde{X} \longrightarrow  {\Bbb C}^{2}  \]
was the pullback of $a : X \longrightarrow  A$ to the universal covers of $X$
and $A$. Consider now a point $x$ in $N$. It follows from the construction that
$g^{-1}(x)$ is contained  $ g^{-1}(N)$. We conclude that  $l(g^{-1}(x))$ has
discrete image under $s$ and therefore it has a discrete image under the map
\[\widetilde{a}:\widetilde{X} \longrightarrow \widetilde{J}.\]
In this way we obtain a fiber $V=\widetilde{a}(g^{-1}(x))$ of the map $v$ such
that $\dim_{\Bbb C}V=0$. Therefore $\dim_{\Bbb C} J \leq 2$ and we are in case
b) of the theorem.
\end{enumerate}

So far we have shown that  the map
\[ g: \widetilde{X} \longrightarrow \Bbb C^{2}\]
has no multiple fibers over an open set  whose complement is of codimension
greater than two. Now we show that in this case part a) of the theorem holds.
To finish the argument we use the discussion of Gromov and Schoen \cite{1}
paragraph 9 about the  Stein factorization for nonproper varieties. The only
thing we need to check is if the coresponding leave space is Hausdorff.  This
follows easily from Lemma 9.3 in \cite{1}. After that the   version of the
Stein  factorization theorem from (see Theorem 3  \cite{33}) implies that the
map $g$ exists, namely we  have $g$ factoring  as
\[ g': \widetilde{X} \longrightarrow \Bbb C^{2}\]
and a finite covering
\[ m: {\Bbb C}^{2} \longrightarrow  {\Bbb C}^{2}.\]
To show that the map $g$ has connected fibers we need to show that   $m$ is an
isomorphism. But since the branch locus of $m$ produces a multiple fibers of
$g$ we conclude that $m$ is etale outside some subset in $\Bbb{C}^{2}$ of
codimension two. Then due  to the purity theorem ( SGA) we conclude  that $m$
is etale. But $\pi_{1}({\Bbb{C}}^{2} ) = 0$ and hence $m$ is biholomorphic.

 Another way to see that $m$ is biholomorphic is to observe that since $Q$ has
real codimension at least 4 then ${\Bbb C}^{2} \setminus Q$ is simply connected
and therefore  it has one, namely the trivial nonramified covering. Therefore
$g$ has connected fibers in codimension two. We finish the proof by observing
that semicontinuity ensures  all fibers of $g$ are connected.
\hfill $\Box$

Now we state the Clemens-Lefschetz-Simpson theorem for   arbitrary $n$. We
sketch a proof emphasizing the details  where it differs  from the case $n=1$.

Recall that $\alpha _{1}, \alpha _{2}, \ldots ,\alpha _{n}$ are $n$ are
holomorphic one forms on $X$ and the map

\[g: \widetilde{X} \rightarrow \Bbb C^{n} \]

was defined on $\widetilde{X}$ by integrating  them  on $ \widetilde{X}$.

\begin{theo}[Clemens-Lefschetz-Simpson] Let $\alpha _{1}, \alpha _{2},\ldots
,\alpha _{n}$ be linearly independent  holomorphic one forms on $X$.  Then one
of the following cases holds:
\begin{list}{\alph{inner})}{\usecounter{inner}}
\item The map $g$ has connected fibers.
\item  There is a holomorphic map with connected fibers from $X$ to the
projective normal variety $Y$ such that $ 1 {\leq }\dim Y {\leq n} $ and the
forms $\alpha _{1}, \alpha _{2}, \ldots ,\alpha _{n}$  are pullbacks from $Y$.
\item The third case is a combination of the previous two. Namely  after a
linear change of the coordinates in ${\Bbb C}^{n}$ some of the forms    $\alpha
_{1}, \ldots ,\alpha _{k}$ come as a pullbacks of a map  $X \longrightarrow Y$,
where $Y$ is an algebraic  variety  $1  \leq  \dim Y  \leq k < n$. The  rest of
the forms  give us a map
$g': \widetilde{X} \rightarrow {\Bbb C}^{n-k} $, with connected fibers.
\end{list}
\end{theo}
{\bf Proof.} The proof splits again into two cases.

\begin{enumerate}
\item  $Q$ contains a hyperplane - then we reduce to the Simpsons's version of
the theorem.
\item  $Q$ does not contain a hyperplane. We argue in the same way as in the
case of $n=2$ . Instead of applying the Castelnuovo De Francis theorem we apply
a generalization of it due to Z. Ran \cite{R}.

\begin{theo}[Z. Ran] Let $\alpha _{1}, \alpha _{2},\ldots ,\alpha _{n}$ be
linearly independent  holomorphic one forms on $X$ and such that
\[ \alpha _{1} \wedge \alpha _{2} \wedge \ldots \wedge \alpha _{n} =0.\]
Then there exists a complex torus $A$ and an analytic map
\[f : X \longrightarrow A ,\]
such that $f(X)$ is a proper annalytic subvariety of $A$ , $\dim f(X)\leq n$,
and the forms \linebreak $\alpha _{1}, \alpha _{2}, \ldots ,\alpha _{n}$ are
pullbacks from $A$.
\end{theo}
\end{enumerate}
\begin{flushright}
$\Box$
\end{flushright}

In  section 5 we give a version of the Clemens-Lefschetz-Simpson
theorem for quasiprojective varieties.

\section{Harmonic maps to  buildings-basic facts}

\subsection{}

Following Gromov-Schoen   and Korevaar-Schoen-Jost  we  briefly describe some
facts of the theory of the harmonic maps to buildings, which are  essential for
our discussion. For a  more detailed treatment one can look at their original
papers \cite{1}, \cite{16}. Also in  this section we define  a new object:  a
spectral covering, which we assign to every nonconstant harmonic map to a
building.

Suppose from now on that $X$ is a smooth projective variety and $\widetilde{X}$
is the universal cover of $X$. Let $G$ be a simple Lie group of rank $r$   over
${\Bbb Q}_{p}$ and  let $B$ be the corresponding Euclidean building. (For
definitions and more detailed account on  buildings  see \cite{6}.) $G$  acts
on $B$ by automorphisms  and if there exists   a representation of $\pi_{1}(X)$
in $G$, so does  $G$.

Following \cite{1}  one  can define an  energy functional  for every
equivariant continuous nonconstant map $U:\widetilde{X}{\longrightarrow}B$ and
we can ask if there exists  a minimum for this functional - the harmonic map.
The following remarkable result gives the answer to  this question.

\begin {theo}[Gromov-Schoen] If  the action of  $\pi_{1}(X)$ on the  building
has no fixed  point, then there exists an equivariant harmonic map $U :
\widetilde{X} \longrightarrow B$.
\end{theo}

{\bf Proof.} See  \cite{1}. \hfill $\Box$

Another  important result from  \cite{1}  which we will need   is the
following:

\begin {theo}[Gromov-Schoen] The singular set of the equivariant harmonic map
$U : \widetilde{X} \longrightarrow B$ is of   codimension  at least two.
\end {theo}

All these results  hold in the case when the building is locally compact. For
some of our considerations we need the same results in the nonlocally compact
situation. In this case we  use  the following:

\begin {theo}[Korevaar -Schoen] The previous two theorems hold in the case of
length spaces of non-positive curvature, including nonlocally compact
buildings.
\end{theo}
{\bf Proof.} The theorem follows from proposition 2.6.5 in  \cite{16}. The
actual proof was given by R. Schoen in his course at IAS Princeton 1992-1993
and will appear elsewhere. R. Schoen informed me  that the same result was
proven J. Jost. \hfill $\Box$

Using this map $U$ we  construct a new object - a spectral covering. As we will
see later this construction is not unique.

 In every apartment $Ap$ of the building we have  locally a canonical choice of
coordinates modulo similarity maps  and actions of the affine Weyl group
$\widetilde{W}$. This means that locally we can choose  one coordinate $z_{1}$
 and consider also $e_{1}(z_{1}) , \ldots ,e_{w}(z_{1})$, where $e_{1}\ldots
e_{w}$ are the elements of the usual Weyl group, and $w$= number of the
elements in $W$.
 Let  $U^{*} T^{*}Ap $ be the pullback of the cotangent bundle of $Ap$  and let
$d$ denote the corresponding exterior derivative in $T^{\star}Ap $. Then we can
define, modulo the action of $W$,  the differential  forms $de_{1}(z_{1}),
\ldots ,de_{w}(z_{1})$ globally on the whole building. Observe that after
taking  differentials we can make everything invariant under the translations
of   the affine Weyl group. Consider  the complexified  differentials
 $U^{\star}de_{1}(z_{1}) , \ldots ,U^{\star}de_{w}(z_{1})$. Due to the
harmonicity of $U$ the (1,0) part of each of the forms $U^{\star}de_{1}(z_{1})
, \ldots, U^{\star}de_{w}(z_{1})$ is
holomorphic. (This follows essentially from the fact that  due to  the Corlette
vanishing theorem (see \cite{1}), $U$  is actually pluriharmonic outside the
singular set). We will denote these  new  holomorphic one forms  again by
$U^{\star}de_{1}(z_{1}) ,\ldots ,U^{\star}de_{w}(z_{1})$. Since $U$ is an
equivariant map these holomorphic one forms descend to    forms on $X$, defined
modulo a $W$-action. Let $h_{1}\ldots h_{r}$ be  a basis  for  all invariant
polynomials of $G$. We  apply  them to $U^{\star}de_{1}(z_{1}),\ldots ,
U^{\star}de_{w}(z_{1})$.  According to \cite{1} the map $U$ is Lipschitz  and
therefore we get that the  holomorphic differentials     $h_{1}, \ldots ,
h_{r}$ are bounded. Using  Theorem 3.2  we conclude that  $h_{1}, \ldots ,
h_{r}$ extend to a holomorphic differentials on the whole $X$, namely $h_{1},
\ldots , h_{r}$ extend to elements of   $H^{0}(X,{\rm
Symm}^{d_{1}}\Omega^{1}_{X}), \ldots, H^{0}(X,{\rm
Symm}^{d_{r}}\Omega^{1}_{X})$ respectively.

\begin{rem} {\rm If the group $G$ is not simply connected one needs include in
$W$ the fundamental group of $G$ as well.}
\end{rem}

\subsection{}

Here  we give one way of defining a spectral covering:

We define a  spectral covering corresponding to the  elements  $h_{1},\ldots
,h_{r}$ from $H^{0}(X,{\rm Symm}^{d_{1}}{\Omega^{1}_{X}}),$ $\ldots,
H^{0}(X,{\rm Symm}^{d_{r}}{\Omega^{1}_{X}})$ respectively  to be  the zero
scheme ${\frak S}$ of the section $\det{(\rho(\theta) - \lambda\cdot id_{\bf
V})}$ , where $id_{\bf V}$ is the identity in a vector space $V$ such that
$\dim {\bf V} = w$. Here by $\det{(\rho(\theta) - \lambda\cdot id_{\bf V})}$ we
mean the following :

\begin{enumerate}
\item  For every element of $G$ the coefficients of the characteristic
polynomial $\det{(\rho(\theta) - \lambda\cdot id_{\bf V})}$ (here $\lambda $
is a number ) are combinations of the invariant polynomials $h_{1},\ldots
,h_{r}$.

\item  Now if we interpret $h_{1},\ldots ,h_{r}$ as elements from
$H^{0}(X,{\rm Symm}^{d_{1}}{\Omega^{1}_{X}}), \ldots , H^{0}(X,{\rm
Symm}^{d_{r}}{\Omega^{1}_{X}})$ and $ \lambda$  as  the tautological section of
 $T^{*}X$ we obtain $\det{(\rho(\theta) - \lambda\cdot id_{\bf V})}$ as an
element from $H^{0}(T^{*}X,{\rm Symm}^{w}{\pi^{*}\Omega^{1}_{X}})$.
 \end{enumerate}

For more detailed account on spectral coverings  see \cite{8}.

Some of these coverings may be nonreduced  and reducible. In such  cases we
work with $S' = {\frak S}_{\rm red}$. Finally we take $S$ - the equivariant
desingularization of $S'$. (Note that the existence of this desingularisation
follows from the equivariant version of the famous Hironaka's theorem.)

\begin{defi} $S$ is  called the factorizing spectral covering corresponding to
$U$.
\end{defi}

For the proof of the main theorem we  need  only the following three properties
of $S$:

\begin{list}{\arabic{inner})}{\usecounter{inner}}
\item The forms $U^{\star}de_{1}(z_{1}), \ldots , U^{\star}de_{w}(z_{1})$ are
well defined holomorphic one forms on $S$;
\item There is a $W$ action on at least an open set of $S$;
\item $S$ is smooth.
\end{list}

Therefore we can define $S$  as a finite covering of a blow-up of $X$ with the
above three properties. An alternative way of constructing $S$ is  given in
\cite{8} and \cite{9}. We describe this construction now .

First choose $r$  generically linearly independent holomorphic one forms out of
 $U^{\star}de_{1}(z_{1}),\ldots ,$ $U^{\star}de_{w}(z_{1})$. The multivalued
form $U^{\star}de_{1}(z_{1}),\ldots , U^{\star}de_{r}(z_{1})$ defines
an $r$-fold covering  $S'$ in $T^{*}X$. Now we define $S$ to  be the Galois
closure of the function field extension ${\Bbb C}(S')/{\Bbb C}(X)$.

It is clear that   from the above definitions that the spectral covering is
not unique. In the proof of the main theorem we show that the factorizing
properties of a  given representation do not depend on the spectral  covering
we have chosen.

The following fact is going to play an essential role in the proof of the
Shafarevich-Koll\'{a}r Conjecture.

\begin{lemma} The singular set $Q$ of the map $U^{'} : S \longrightarrow B$ is
contained in the union of the zeros of $U^{\star}de_{1}(z_{1}),\ldots ,
U^{\star}de_{w}(z_{1})$.
\end{lemma}
{\bf Proof.} The proof follows from the fact that according to Gromov and
Schoen ${\rm codim}_{S}Q \leq 2$ and that $U^{'}(Q)$ is contained in the faces
of the chambers of $B$. The existence of the intrinsic derivative, proved in
\cite{1} and \cite{16}, implies existence of a kernel $T^{*}S$ and therefore
$U^{\star}de_{1}(z_{1}),\ldots , U^{\star}de_{w}(z_{1})$ are  equal to zero on
$Q$. \hfill $\Box$

\begin{rem} Observe that over some open set in $S$ , namely outside $Q$ we can
use the Castelnouvo-de Francis theorem to get a factorization  (see \cite{2}).

\end{rem}

\subsection{Example}

 Let us consider the case of the group $SL(3,{\Bbb  Q_{p}})$. The Weyl group in
this case is going to be $S_{3}$. Then the
coordinates $e_{1}(z_{1})\ldots e_{6}(z_{1})$  can be thought as
 six roots of unity on the face of $B$, which in this case is nothing else but
an equilaterial triangle.
The complexified  differentials $U^{\star}de_{1}(z_{1}),\ldots
,U^{\star}de_{6}(z_{1})$ are
$U^{\star}de_{1}(z_{1}), \ldots,  U^{\star}de_{3}(z_{1}) ,
-U^{\star}de_{1}(z_{1}), \ldots,  -U^{\star}de_{3}(z_{1})$.

In this case we have

$h_{1}=-(U^{\star}de_{1}(z_{1})^{2}+U^{\star}de_{2}(z_{2})^{2}+U^{\star}de_{3}(z_{3})^{2})$

$h_{2}=U^{\star}de_{1}(z_{1})^{2}.U^{\star}de_{2}(z_{2})^{2}.U^{\star}de_{3}(z_{3})^{2}$.

The spectral  covering $\widetilde{S}$ in this case is  given by the following
section in $H^{0}(X,{\rm Symm}^{6}\Omega^{1}_{X})$:

$\lambda^{6}-h_{1}\lambda^{4}+h_{1}^{2}/2{\lambda^{2}}-h_{2}^{2}$,

where  $\lambda$ is the tautological one form.

If, for example,
 $U^{\star}de_{1}(z_{1})=U^{\star}de_{2}(z_{2})\neq 0$, then the factorizing
spectral covering $S$ is the desingularization of the zeros of the section
$\lambda^{2}-U^{\star}de_{1}(z_{1})^{2}$, which belongs to  $H^{0}(X,{\rm
Symm}^{2}\Omega^{1}_{X})$.

 \section { The main theorem and its proof}
\subsection{}

The main idea is to use the abundance of  holomorphic one forms over the
factorizing  spectral covering and to apply to it the Clemens-Lefschetz-Simpson
theorem we have proven above. This gives the factorization over an open set of
$X$. To  finish the proof we use a new technique developed by J\'anos Koll\'ar
in \cite{30}  and  \cite{31}.

Recall some notation from section 2. Let  $X$ be  a smooth projective variety
and let $\widetilde{X}$ be  the universal cover of $X$. Let $G$ be a simple Lie
group of rank $r$ over an arbitrary local field of   characteristic zero,  and
let  $B$ be  the corresponding Euclidean building ($B$ could be  non locally
finite). $G$  acts on $B$ by isometries  and if we have a representation
$\varrho:\pi_{1}(X) \longrightarrow G $, so does $\pi_{1}(X)$. If  $\pi_{1}(X)$
acts on $B$ without a fixed point  then it follows from the theory of
Gromov-Korevaar-Schoen-Jost  that there exists a nonconstant   equivariant
harmonic map $U:\widetilde{X}\longrightarrow B$.

\begin {theo} Let  $U:\widetilde{X}\longrightarrow B$ be an  equivariant
harmonic map which corresponds to the Zariski dense representation   $\varrho:
\pi_{1}(X) \longrightarrow G$. Then there exists a holomorphic map $x:S
\longrightarrow Y $ where $Y$ is a normal projective variety  such that $1 \leq
 dim(Y) \leq r$. Moreover there exist   $S^{0}$ and  $Y^{0}$ - Zariski open
sets in $S$ and $Y$ respectively such that  $\varrho: \pi_{1}(S^{0} )
\longrightarrow G$ factors through a representation $\varrho{\prime}:
\pi_{1}(Y^{0}))\longrightarrow G$.
\end {theo}

The fact that  $\pi_{1}(X)$ acts on $B$ without fixed points  allows us to
define spectral covering  $\widetilde{S}$ the corresponding to $\varrho:
\pi_{1}(X) \longrightarrow G$. Take the corresponding $S$. Since the map
$U:\widetilde{X}\longrightarrow B$ is a nonconstant map we have  at least one
holomorphic one form  $\alpha _{1} $ on $S$.

 If we have the forms   $\alpha _{1},\ldots ,\alpha _{q}$, $1\leq q \leq r$,
then we are in position to apply the  Clemens-Lefschetz -Simpson theorem.

 Suppose that part  a)  or a) of this theorem holds in our situation.

Let  $\bar{S}$  be  the universal cover of $S$. Consider the map
$g: \bar{S} \longrightarrow {\Bbb C}^{q}$  defined  as in   section 2  by
integration   of  $\alpha _{1}, \ldots ,\alpha _{q}$ over
$\bar{S}$. Let $re(g)$ be the real part of this map. Obviously
$U : \widetilde{X}\longrightarrow B$ provides  us  with the
$\pi_{1}(S)$-equivariant harmonic  map $U' :\bar{S}\longrightarrow B$.

\begin{lemma} The fiber of $re(g) :\bar{S} \longrightarrow {\Bbb R}^{q}$ maps
to a point under $U'$.
\end {lemma}
{\bf Proof.} We  use an argument of  \cite{1}, namely the local version of the
Stein factorization theorem, which says that  in a small neighborhood $\Omega$
in  $\bar {S}$ the map $g$ decomposes into a composition of a holomorphic map
$\Omega\longrightarrow D$, where $D$ is an  open ball in  ${\Bbb C}^{q}$,
followed by a harmonic map $u : D \longrightarrow B$. So we get that each
intersection of the fiber of $re(g)$ with  a small neighborhood $\Omega$ goes
to a point  in $B$. But now using the fact that the fiber is connected and
all the maps are continuous we see that the lemma holds. Observe  that this
argument also  works in the neighborhoods around  the critical points of this
map. What one does is first approximate the singular fibers by nonsingular
ones, and then  use the fact that $g$ is a continuous map. \hfill $\Box$

Now using the  above Lemma  we  construct an isometry ${\Bbb R}^{q}
\longrightarrow   B$. Before we do  this we  show how we  are going to use this
isometry.

Recall that we are in case a) or c)  of the
Clemens-Lefschetz-Simpson theorem with $l=q$. The image of
$U$ is all of ${\Bbb R}^{q}$. But this means that the  action  of  $\pi_{1}(X)$
on the building  $B$   preserves ${\Bbb R}^{q}$, so the action on the building
at infinity  preserves a flat.  Hence  it is contained  in some parabolic
subgroup in $G$. This  contradicts the Zariski density  of $\varrho$.  In fact
one can prove this without going to the building at infinity. Fixing    ${\Bbb
R}^{q}$  is equivalent to fixing ${\Bbb R}^{q}$ pointwise  up to  the action of
some finite group and this is again a  contradiction.

This rules out  the cases a) and c)  of the
Clemens-Lefschetz-Simpson theorem.

Now we prove the existence of the isometry ${\Bbb R}^{q}\longrightarrow B$.

\begin{lemma} In the  cases when  part a) or c)  of the
Clemens-Lefschetz-Simpson theorem holds we always have an isometry ${\Bbb
R}^{q} \longrightarrow B$.
\end {lemma}
{\bf Proof.} The case $q=1$  was proven by Simpson \cite{4}. Following our
proof of the
Clemens -Lefschetz -Simpson theorem we are going to work for simplicity with
the case $q=2$.

Start with the map $f:\bar{S}\longrightarrow{\Bbb C}$ defined by integration of
$\alpha _{1}$ over $\bar{S}$. Since the singular sets of the map $f$ are
compact sets  in  $\bar{S}$ and the fibers $F$ of $f$ are connected (recall we
are in case a) or c)  of the
Clemens-Lefschetz-Simpson theorem), we can find a real line ${\Bbb R}$ in
${\Bbb C}$  over which the map $f$ has only smooth points. Now through  every
point of this  real line  we can define another real line  in the fiber of the
map $t : F \longrightarrow {\Bbb C}$. Here we again use  the fact that  the
singular sets of the map $t$ are  compact sets  in  $F$ and the fibers $F'$ of
$t$ are connected. This means that we can not only find such a line
 but also require that it  passes through the point of our initial line ${\Bbb
R}$ in ${\Bbb C}$. We can do this in a continuous way  using the fact that the
singular sets of the maps above are in codimension at least two. Since every
${\Bbb R}$-bundle over ${\Bbb R}$  is trivial we have a smooth ${\Bbb R}^{2}$
in $\bar{S}$.
 Now using Lemma 4.1 we find smooth map $\Bbb R^{2}\longrightarrow B$. To show
that this is an isometry we  apply  Simpson's argument, which uses that the
forms $\alpha _{1}$  and $\alpha _{2}$ are defined by $U$. Thus  the
differentials of $U$ and $re(g)$ are the same and therefore  the differential
of $U$ is equal to the identity.
$\Box$

 From now on we are going to work only with case B) of the Clemens -Lefschetz
-Simpson theorem. We need  first to construct the factorization  map. Part b)
of the Clemens-Lefschetz-Simpson theorem gives  us a map  $S\longrightarrow Y$.
 So we have a $\pi_{1}(S)$ - equivariant harmonic map
$U\prime:\bar{S}\longrightarrow B$. Let $\widetilde{Y}$  be the universal cover
of  $Y$.

\begin{lemma} The map  $U\prime:\bar{S}\longrightarrow B$  factors through a
map
$u_{0}:\widetilde{Y}\longrightarrow B$, namely $U\prime=u_{0}.a$, where $a $ is
the map $a:\bar{S}\longrightarrow\widetilde{Y}$.
\end{lemma}
{\bf Proof.}
For simplicity we again consider  the case  $dim_{\Bbb R}  B=2$ and two
holomorphic one-forms  $\alpha _{1}$, $\alpha _{2}$ such that $\alpha
_{1}\wedge \alpha _{2}\neq0$ generically on $\bar{S}$. The forms  $\alpha
_{1}$, $\alpha _{2}$  come from holomorphic  one-forms $\beta _{1}$, $\beta
_{2}$ on $\widetilde{Y}$. This way  we get that the  map
$r_{1}:\bar{S}\longrightarrow \Bbb R^{2}$, given by integration of $\alpha
_{1}$, $\alpha _{2}$ over $\bar{S}$, factors through the map $r_{2}:Y
\longrightarrow \Bbb R^{2}$, given by integration of $\beta _{1}$, $\beta _{2}$
on $\widetilde{Y}$. Namely we have $r_{1}=r_{2}.a$. But then the fiber of the
map $a:\bar{S}\longrightarrow\widetilde{Y}$ is contained in the fiber $J$  of
the map $r_{1}:\bar{S}\longrightarrow \Bbb R^{2}$. By  Lemma 4.1 $J$ is mapped
to a point under $U$ and so is the fiber of the map $a$. This proves the
factorization.

 $\Box$

Define $X^{0}$  to be the   Zariski dense set in $X$ obtained by throwing away
the branch locus  of the map $S \longrightarrow X $ and the images  in $X$ of
the exceptional sets  in $S$. Let $S^{0}$ be the preimage of $X^{0}$ in $S$.
The fundamental group $\pi_{1}(S^{0})$   maps to a group with finite index in
$\pi_{1}(X^{0})$. The representation $\varrho(\pi_{1}(X))$ is Zariski dense in
 $G$ as is $\varrho (\pi_{1}(X^{0}))$ since $\pi_{1}(X^{0})$ surjects to
$\pi_{1}(X)$. But since we know that $\pi_{1}(S^{0})$  maps to a group with
finite index in  $\pi_{1}(X^{0})$,  we can conclude that
$\varrho(\pi_{1}(S^{0}))$ is Zariski dense in   $G$.

 From the previous lemma we have that the harmonic map
$\bar{U}:\widetilde{S^{0}}\longrightarrow B$ factors through the  map
$\bar{u_{0}}:\widetilde{Y^{0}} \longrightarrow B$, namely
$\bar{U}=\bar{u_{0}}.a$. (Here  $\widetilde{S^{0}}$ and  $\widetilde{Y^{0}}$,
are the universal covers of  $S^{0}$ and $Y^{0}$ respectively.) Using  the
properties of the pullback map for the inclusion $ i: \widetilde{S^{0}}
\longrightarrow \bar{S }$ and  the fact that $\pi_{1}(S^{0})$ surjects onto
$\pi_{1}(S)$ we conclude that the map $\bar{U}:\widetilde{S^{0}}\longrightarrow
B$ is equivariant with respect to  $\pi_{1}(S^{0})$.

   We  show that the action of $\varrho:\pi_{1}(S^{0})$ factors through an
action
of  $\varrho{\prime}:(\pi_{1}(Y^{0}))$. Take $\gamma$ an element of
$\pi_{1}(S^{0})$  such that
  $a_{\star}(\gamma)=1$. Then for any $x$ in  $\widetilde{S^{0}}$ we have
$a_{\star}(\gamma.x)=
  a_{\star}(x)$. We use the same notation: $a$ for both maps
$a:S\longrightarrow Y$ and $a:\bar{S}\longrightarrow\widetilde{Y}$ .
  Then using the equivariance of $U$    we have

\[\varrho(\gamma)\bar{U}(x)=\bar{U}(\gamma.x)=\bar{u_{0}}.a(\gamma.x)=\bar{u_{0}}.a(x)=\bar{U}(x).\]

	In the same way one can see that  $\bar{u_{0}}$ is equivariant for the action
of $\varrho{\prime}(\pi_{1}(Y^{0}))$ on $B$.

 \[\bar{u_{0}}(a_{\star}(\gamma)a(x))=\bar{u_{0}}a(\gamma x)=\bar{U}(\gamma
x)=\varrho(\gamma)\bar{U}( x)=
  \varrho(a_{\star}\gamma)\bar{U}(ax).\]

  We need to show that  $\varrho(\gamma)=1$ if $\gamma$ is in $Ker (a
_{\star})$. Recall that $a$ has connected fibers, so the map
$a_{\star}:\pi_{1}(S^{0})\longrightarrow \pi_{1}(Y^{0})$ is a surjective map.
But $Ker (a _{\star})$ is a normal subgroup in $\pi_{1}(S^{0})$ thus
$\varrho(Ker (a _{\star}))$ is normal  in $G$.  Since $G$ is a simple Lie group
over $K$ there are  two  possibilities:

1) $\varrho(Ker (a _{\star}))=G$;

2) $\varrho(Ker (a _{\star}))$ is contained in the center $Z(G)$.

 Suppose that we are in case 1). But from the computation above we know that
$\varrho(Ker (a _{\star}))$ fixes $\bar{U}(X^{0})$, an open  set in $U(X)$,
which is a contradiction since $G$  does not fix  any point in $B$.

In case 2) we obtain that $\varrho(\gamma)=1$ since we are working only with
centerless representations.

$\Box$

Let $\varrho:\pi_{1}(X) \longrightarrow G$ be a  Zariski  dense nonrigid
representation of the   fundamental group of $X$ to some  complex  simple Lie
group $G$.  Then there exists a  curve $I$ in the moduli space of
representations passing through
  $\varrho:\pi_{1}(X) \longrightarrow G$ in direction of which $\rho$ deforms
nontrivialy. Since the moduli space of representations is an affine variety the
curve $I$  is affine as well and let $\bar{I}$ be the compactification of $I$.
Consider now the representation  $\bar{\rho}$ to $\bar{G}$, where $ \bar{G} $
is defined over the field of  fractions of  $I$. Let $O_{p}$ be the
localization at some point $p \in \bar{I} \ I$. Let $Z[T]$ be an extension of
the ring of integers of  $O_{p}$ which contains all coefficients of
$\bar{\rho}$ localized at the point $p$. Consider now $I$ as a curve over $Spec
Z[T]$. The fact that $\varrho:\pi_{1}(X) \longrightarrow G$ is a   Zariski
dense nonrigid representation  implies that $\bar{\rho}$ is   a   Zariski
dense nonrigid representation. Observe that there  are still going to be a
point $q$ in the  curve $\bar{I}$ which is not in $Spec Z[T]$. Let
$\chi(\bar{\rho})$ is a character of  $\bar{\rho}$. Moving over $\bar{I}
\subset Spec Z[T]$ we see that $\chi(\bar{\rho})$ is unbounded at $q$.
Therefore we conclude that  the  representation $\bar{\rho}$ to $\bar{G}$ is
not contained in any bounded subgroup in $\bar{G}$.

The following corollary is an easy consequence  of the previous result:

  \begin{corr}  Let $U$  be  an  equivariant harmonic map
$U:\widetilde{X}\longrightarrow B$ and   $\varrho:\pi_{1}(X) \longrightarrow G$
be a Zariski dense nonrigid representation, and  let $X$ be  a projective
variety. Then one of the following possibilities holds:

  A) The $Im(\varrho(\pi_{1}(X))$ is contained in a maximal compact subgroup
in $G$ or

B) There exists a holomorphic  map $S \longrightarrow Y$ and a map
$u_{0}:S^{0} \longrightarrow Y^{0} $ with $S^{0}$ and $Y^{0} $ Zariski open
sets in  $S$ and $Y$ respectively,  such that $\varrho( \pi_{1}(S^{0})) $
factors through a representation $\varrho{\prime}: \pi_{1}(Y^{0})
\longrightarrow G$, where $Y$ is a normal projective variety  such that $1 \leq
 dim(Y) \leq r$.

 \end{corr}

Using this corollary we give a proof of Theorem 1.2. It follows from remark 3.1
 and \cite{2} , paragraph 4 that using the properties of the spectral covering
we can mod out  $S^{0}$ by  $W$ and  get a map

\[h:X^{0}\longrightarrow Y^{0} ,\]

, which comes from the morphism $h:X \longrightarrow Y/W $, where  $X^{0}$ is
an open set  in $X$ and $Y^{0}$   is an open set in $Y/W$. From now on we
denote all modifications of $Y$, namely blow ups and finite nonramified covers,
by $Y$ if not stated otherwise.

Denote the generic fiber of $h:X^{0} \longrightarrow Y^{0}$  by $Z^{1}$.

 Let us first  resolve the singularities of $Y$. This might  change the
fundamental group of $Y$  but we need only that  the fundamental group of $X$
does not change, which  follows from the fact that   $X$ is a smooth projective
variety.

We finish the proof of the main theorem by generalizing 4.8.1. in \cite{30}.
According to corollary 4.1 we have the following sequence:

\[ \pi_{1}(X^{0}) \longrightarrow \pi_{1}(Y^{0}) \longrightarrow G.\]

What we need to show is that $\pi_{1}(Z^{1})$ belongs to the kernel of the map

\[ \pi_{1}(X^{0}) \longrightarrow \pi_{1}(Y^{0}).\]

We cannot do that directly on $X$ but by generalizing 4.8.1. in \cite{30} we
show that this is possible on some finite nonramified covering of $X$, $X(T)$.

Observe that

\[ im[\pi_{1}(X^{0}) \longrightarrow \pi_{1}(Y^{0})]\subset Ker \rho.\]

\begin{defi} Let $P$ be  a subgroup of $ \pi_{1}(X)$ and let $X(P)$ be a
covering of $X$ with fundamental group $P$. Consider now the Stein
factorization of the map $X(P) \longrightarrow Y$ and
define $X(P) \longrightarrow Y(P) $ to be the map with connected fibers in this
factorization.

\end{defi}

 Note that  $Y(P)$ will be an analytic variety even when  the covering $
Y(P)\longrightarrow Y$  is  infinite. This follows from the most general
version of the Stein  factorization theorem (see Theorem 3  \cite{33}).

\begin{defi} Define $ \Omega$ to be the intersection of all subgroups $P$ in $
\pi_{1}(X)$ such that $H\subseteq P $, where
$H=im[\pi_{1}(Z^{1})\longrightarrow \pi_{1}(X)]$ and the covering  $
Y(P)\longrightarrow Y$ has finite ramification indexes.

\end{defi}

Such an $\Omega$  is well defined due to the fact  that there exists at least
one such a $P$ - $\pi_{1}(X)$.

Define $K=Ker (\varrho) $.

\begin{lemma} The covering $Y(K) \longrightarrow Y$  has finite ramification
indexes.

\end{lemma}

{\bf Proof.}

After  intersecting  $X$  with sufficiently many generic hyperplanes   we get a
 finite ramified covering $X \cap H \longrightarrow Y$. This covering obviously
has finite  ramification index. But since the covering  $X(K)\longrightarrow X$
is nonramified  we know the indeces of the covering $Y(K) \longrightarrow Y$
are  also finite.

$\Box$

\begin{lemma} There exists a finite index subgroup $T$ in $\pi_{1}(X)$ such
that:

\[H  \subseteq  R  \subseteq  \Omega \subseteq \pi_{1}(X ),\]

where $R$  is  the kernel of the map $\pi_{1}(X(T))\longrightarrow
\pi_{1}(Y(T))$.

\end{lemma}

{\bf Proof.} The proof is the same as  in 4.8.1 \cite{30}.

$\Box$

To finish  the proof of 1.2  we  need only   to observe that

\[ R \subseteq \Omega \subseteq K \]

and this gives us the complete factorization

\[X(T) \longrightarrow Y(T).\]

As an almost immediate corollary we have:

\begin{corr}

 Let $\varrho:\pi_{1}(X) \longrightarrow G$ be a  Zariski  dense nonrigid
representation of the   fundamental group of $X$ to some  complex  simple Lie
group $G$. Then  there exist:

1)  a finite etale cover $X'$ of a blow up of $X$;

2) a smooth projective variety $Y$ of positive dimension    $l$ less then or
equal to the rank $r$  of $G$ over $\Bbb{C}$;

3) and a  holomorphic map $h:X' \longrightarrow Y $

such that $\varrho:\pi_{1}(X') \longrightarrow G$ factors through a
representation of $\pi_{1}(Y) $.

\end{corr}

{\bf Proof.}   The fact that $\varrho:\pi_{1}(X) \longrightarrow G$ is  a
Zariski  dense nonrigid representation of the   fundamental group of $X$ to
some  complex  simple Lie group $G$ implies that there exists a curve in the
moduli space of representations passing through
  $\varrho: \pi_{1}(X) \longrightarrow G$, which  intersects  infinity in the
moduli space of representations in the point $p$. Let $O_{p}$ be the local ring
of this point as a point of the curve described above. This way we obtain  a
representation $\bar{\rho}$ to $\bar{G}$, where $ \bar{G} $ is  defined over
the  field of fractions of the completion of $O_{p}$, which is also  Zariski
dense and is not contained in any bounded subgroup in $\bar{G}$.  We finish the
proof  by applying  Theorem 1.2 .

$\Box$

\begin{rem} In the last section we give a different proof of the
last part of  Theorem 1.2. There we  use  the theory of the Shafarevich maps.

\end{rem}

\begin{rem} Observe  that if we work with group of finite index in
$\pi_{1}(X)$ we get again the  conclusions of the main theorem.
\end{rem}

\section{The quasiprojective case.}

In this section we are going to work with $X=X_{1}\setminus D$ - a
quasiprojective variety , where $D$ is the divisor at infinity and $X_{1}$ is
the compactification of $X$. Our goal is to prove theorem similar to that of
the previous section. Of course to start the whole procedure we need  to show
that there exists an equivariant continuous map of finite energy
$U:\widetilde{X}\longrightarrow B$. To be able to show this we  require that
our divisor at infinity is a divisor with normal crossing with unipotent
monodromy around it.  We can cover $X$ by finitely many open sets since it is
compactly embedded. Using the fact that we have a unipotent monodromy at
infinity we see that the fundamental groups of the  open sets  which cover the
divisor at infinity  are contained in  maximal compact subgroups.  A simple
computation shows that by simultaneous conjugation  by  some elements we can
make all of the finitely many generators of these groups have  integer
coefficients
                             and hence they are contained in the same  maximal
compact subgroup. Using the fact that the maximal compact subgroups fix a point
in $B$  we get an  equivariant continuous map of finite energy from every open
set which covers the divisor at infinity to $B$, namely the  map to the fixed
point.  Using the standard center of  mass construction we obtain the
equivariant continuous map of finite energy $U:\widetilde{X}\longrightarrow B$.

 We are going to make sense of spectral covering in the quasiprojective case.
First we formulate and give a sketch of the proof of  Corlette vanishing
theorem in the quasiprojective case.

Let $X$ be a quasiprojective variety with universal covering  $\widetilde{X}$.
Consider now a representation of the fundamental group of $X$ to some Lie group
defined over an arbitrary nonarchimedian field.  Let $B$ be  the corresponding
Euclidean building. Let us assume also that the divisor at infinity in $X$ is a
divisor with normal crossings and the monodromy around it is unipotent. We have
the following theorem,  the proof of which is similar  to the proof of the
original theorem of Gromov and Schoen.
It is clear that in  the situation above, provided that
$\pi_{1}(X)$ acts on $B$  without fixed points, to every Zariski dense
representaion we assign  a harmonic map $U:\widetilde{X}{\longrightarrow}B$
for which  the whole theory of Gromov,  Korevaar, Schoen and Jost works.

Now  take a regular point $x_{0}$ of $\widetilde{X}$ (for the definition of
regular point see \cite{1}). Then the image of the ball $B-\sigma(x_{0})$ is
contained in at least one flat $F$ in $B$. Let $\nabla $ denote the pullback
connection and let $d_{\nabla}$ denote the corresponding exterior derivative
operator on p-forms with values in $U^{*}TF$. Let ${\delta} _{\nabla}$ denote
its formal adjoint. The differential $dU$ then defines a 1-form with values in
$U^{*}TF$, and  we have the following:

\begin {theo}(Corlette ) Let $X$ be a quasiprojective variety, let  $ \omega $
be  a parallel p-form on  $\widetilde{X}$, and  let $U$  be a harmonic map
$U:\widetilde{X}{\longrightarrow}B$.
Then in the neighborhood of  $x_{0}$, a regular point for $U$, the form $w
{\wedge}dU$
satisfies ${\delta} _{\nabla}(w {\wedge} dU){\equiv}0$.

\end {theo}

{\bf Proof.} The statement of this theorem is local  so the proof is almost the
same as in \cite{1}. First we exhaust $X$ by compact sets $X_{i}$. After that
one needs to choose the right cutfunction and apply to each  $X_{i}$ theorem
7.2 from \cite{1}. But since the statement is local we just choose the same
functions as in \cite{1}.

$\Box$
\begin{corr} In the situation above $U$ is pluriharmonic.

\end{corr}

{\bf Proof.} Proof is as in \cite{1}.
$\Box$

First we define the spectral covering . As in the compact case we have the one
forms \linebreak $U^{\star}de_{1}(z_{1}),\ldots , U^{\star}de_{w}(z_{1})$
defined on  $X=X_{1}\setminus D$ up to action of the Weyl group. (Here $z_{1}$
is a an arbitrary choice of a coordinate on a given chamber  of  $B$ , which we
extend using the action of the group.) Again we apply to them the basis of the
$G$-invariant for $G$ polynomials  and obtain the forms  $h_{1}, \ldots ,
h_{r}$. Due to the fact that $U$ is a finite energy map it follows from
\cite{1} that the forms $h_{1}, \ldots , h_{r}$ are $L_{2}$ bounded and
therefore have at most log poles. The unipotency of the loops around infinity
gives us that the residues of  $h_{1}, \ldots , h_{r}$ are rational. Therefore
over some finite covering of $X_{1}$ they are holomorphic. Using the procedure
of section 3 we build now the spectral covering $S$ over $X_{1}$.

The theorem above gives us the forms $U^{\star}de_{1}(z_{1}),\ldots,
 U^{\star}de_{w}(z_{1})$ defined on  $X=X_{1}\setminus D$ up to  an action of
the Weyl group. Therefore over the spectral covering $S$ we obtain the forms
$\alpha _{1}, \alpha _{2},...,\alpha _{r}$, where $r$ is the rank of the group
$G$. Observe that the forms $\alpha _{1}, \alpha _{2},...,\alpha _{r}$  are
$L_{2}$ bounded and therefore have at most log poles. This can be seen in the
following way: Using  \cite{1}(Theorem 6.3) we see that  a loop around $D$
fixes a point $s$ in $B$. Therefore  the map $U:D \longrightarrow B$  gives us
a  weakly subharmonic map $d^{2}(U(x),s)$ and then the same  arguments as in
\cite{1}( Theorem 6.3)   imply the above statement.

Now we state the quasiprojective version of  the  Clemens -Lefschetz -Simpson
theorem.

Observe that if $l:S \longrightarrow X$ is the spectral coveing then
$S=S^{1} \setminus l^{*}(D)$ .

 Following Iitaka \cite{I} we introduce the Albanese map $alb:S \longrightarrow
Alb(S)$ for  quasiprojective varieties. As in the case of projective varieties
it is defined by integrals of holomorphic one forms on $S$. Here $Alb(S)$ is a
semiabelian variety- a group extension of  $Alb(S^{1})$ by
$(\Bbb{C}^{*})^{\times l}$.

Define also
\[ A=Alb(S)/B.\]

We have the natural map

\[a:S \longrightarrow  A.\]

Here $B$ is again the maximal abelian subvariety over which  the forms

$\alpha _{1}$, $\alpha _{2}$,...,  $\alpha _{n}$ are zero as a forms on
$Alb(S^{1})$.

Let us also define the map

\[g:\widetilde {S}\longrightarrow \Bbb{C}^{n}\]

as a pullback of $a$. In this situation we have the following :

\begin{theo}(Clemens - Lefschetz - Simpson) Let $S$ be  a smooth
quasiprojective variety. Let  $\alpha _{1}, \alpha _{2},...,\alpha _{r}$ be
holomorphic one forms on $X$  with at most log poles at  $l^{*}(D)$.  Then one
of the following cases holds:

A) The map $g$ has connected fibers.

B) There is a holomorphic map with connected fibers from $S$ to the  projective
normal variety $Y$ such that $ 1 {\leq }dimY {\leq n} $ and the forms $\alpha
_{1}, \alpha _{2},...,\alpha _{n}$  are pullbacks from $Y$.

C)The third case is a combination of the previous two. Namely  after a  linear
change of the coordinates in $\Bbb {C}^{n}$ some of the forms    $\alpha
_{1},...,\alpha _{k}$ come as a pullbacks of a map  $X \longrightarrow Y$,
where $Y$ is an algebraic  variety  $1 {\leq }\dim Y {\leq k}<n $. The  rest of
the forms  give us a map
$g^{'}: \widetilde{X} \rightarrow \Bbb C^{l} $, where $l \leq n-k $,  with
connected fibers.

\end{theo}

The proof is similar to the proof in the projective case.

Now we formulate the version of our main theorem for  a quasiprojective variety
$X$.

\begin{theo} Let $X=X^{1}\setminus D$ be a smooth quasiprojective variety and
$\varrho:\pi_{1}(X) \longrightarrow G$ be a
 Zariski  dense representation of the   fundamental group of $X$ with unipotent
monodromy around $D$, where $D$ is a divisor with a normal crossing. Let $G$ be
a   simple Lie group over $K$.  Then

A. either  the image of $\varrho$ is  in a maximal compact subgroup of $G$ or,

B. there exist:

 1) a blow up   $ X' $ of a finite etale cover of $X^{1}$;

2) a smooth projective variety $Y$ of positive dimension    $l$ less  than or
equal to the rank $r$  of $G$ over $K$;

3) A holomorphic map $h:X' \longrightarrow Y $ such that $\varrho:\pi_{1}(X')
\longrightarrow G$ factors through a representation of $\pi_{1}(Y) $ and such
that the pullback of $D$ in  $X'$ is a pullback from a divisor on $Y$.

\end{theo}

We give an application of the above theorem,  which was suggested by J.
Koll\'ar.

According to N.Mok (see \cite{21})  every real   Zariski dense discrete
representation  in $SL(n,\Bbb{C})$ of noncompact type of the fundamental group
of any   compact K\"{a}hler manifold after some blow up and  finite nonramified
covering factors though the representation of the fundamental group of
projective algebraic variety of general type.

Let   $X=X^{1}\setminus D$ be a quasiprojective variety such that $X^{1
}$ has  Kodaira dimension zero and let $\varrho:\pi_{1}(X) \longrightarrow
SL(n,\Bbb{C})$ be  a real  Zariski dense discrete representation  of the
fundamental group of
$X$. The hypothesis of the   Mok's theorem send us in case B) of the above
theorem. Since $ X' $ is  a finite etale cover of blow up of $X^{1}$  it has
also Kodaira dimension equal to zero . Mok's theorem also tells us that $Y$
from  part B) of the previous theorem is an  algebraic variety of general type.
Therefore we obtain a holomorphic map $h:X' \longrightarrow Y $ from a variety
with Kodaira dimension zero to a   variety of general type and this impossible,
due to a theorem of Kawamata \cite{27}. We have obtained the following:

\begin{corr}Let  $X=X^{1}\setminus D$ be a quasiprojective variety  such that
$X^{1}$ has  Kodaira characteristic zero. Then:

\[\varrho:\pi_{1}(X) \longrightarrow SL(n,\Bbb{R})\]

is a finite group.
\end{corr}

\section{Factorization theorems for complexes of groups.}

In \cite{1} Gromov and Schoen  proved that if the fundamental group of a
quasiprojective variety $X$  admits a decomposition as an amalgamated product
of groups,  then $X$ admits a surjective holomorphic map to an algebraic curve.
In this section we extend this result to higher dimensional complexes of
groups. This  gives a partial answer to a question  Gromov stated in \cite{22}
(section 7). In what follows we are working with the negatively curved
complexes of groups  defined by Benakli in \cite{32}.

 The key idea in  \cite{1}  is that the  Baas-Serre theory  assigns  to an
amalgamated product of groups  a tree  to which one applies the theory of
harmonic maps. According to \cite{17} and \cite{18}(Theorem 6.4) we can assign
to every negatively curved two dimesional complex of groups a connected
negatively curved simplicial cell  complex of dimension 2  with a finite set of
isometry types of cells $T$. Let   $G(T)$   be the universal covering of $T$.
The fundamental group $\pi_{1}(T)$ (see \cite{18})  acts on  $G(T)$ by
simplicial isometries .  Let  $X$ be  a projective variety such that there
exists  a surjective homomorphism $\theta:\pi_{1}(X)\longrightarrow
\pi_{1}(T)$. Therefore $\pi_{1}(X)$ acts on $G(T)$ without fixed points. Using
the fact that  $G(T)$ is contractible (\cite{17} and \cite{18}) we  apply the
Korevaar-Schoen theorem to   get a harmonic map
$U:\widetilde{X}\longrightarrow G(T)$. In this section we work only with $G(T)$
of the type defined by Benakli in \cite{32} paragraph 7. Namely they are
constructed by using the baricentric subdivision of a regular hyperbolic
polyhedra.

The   following theorem is a consequence of the previous sections:

\begin{theo}  Every representation of the  fundamental  group of projective
varieties onto the fundamental  group of  some negatively curved 2-dimensional
polyhedra, comes from the representation of the fundamental group of an
orbicurve.
\end{theo}

{\bf Proof.} The only thing we need to show is how to construct the spectral
covering.
To do  this we  first  choose  local coordinates on $G(T)$ .  Using the
symmetries, defined on $G(T)$  (see \cite{17} and \cite{18})  we can extend
these coordinates  to the all of  $G(T)$.  To obtain holomorphic forms, and
consequently a spectral covering,  we need  vanishing theorems. But since in
the  case of  $G(T)$ we do not have  a statement  about the codimension of the
singularities of the corresponding harmonic map, we need to do some extra work,
which is the essence of the theorem. First we use the following theorem
(\cite{32}):

\begin{theo} (Benakli) $G(T)$ can be embedded isometrically  in a hyperbolic
space $\Bbb{H}^{3}$
in such a way that $G(T)$ is a "deformation retract " of  $\Bbb{H}^{3}_{1}$.

\end{theo}

 Here we denote  by  $\Bbb{H}^{3}_{1}$ the 3-hyperbolic space $\Bbb{H}^{3}$
with the vertexes  of every ideal polytop  thrown away. Observe that
$\Bbb{H}^{3}_{1}$ is  contractible.

We prove a version of the vanishing theorem for $G(T)$ closely following
\cite{1}. We  omit the details of  the argument in \cite{1}, emphasizing only
the differences.

Let us first mention  that  Benakli's construction  can be made  $\pi_{1}(T)$ -
equivariant.  For our argument we need a special kind of retraction. Namely to
show that the    harmonic maps $U_{t}:\widetilde{X}\longrightarrow
\Bbb{H}^{3}_{t}$ exist for every $t$, we need to make sure   that the curvature
of  $\Bbb{H}^{3}_{t}$ is nonpositive for every $t$.
We do that by retracting  $\Bbb{H}^{3}_{1}$ to $G(T)$ equidistantly. Observe
that the  $\Bbb{H}^{3}_{t}$ are going to be singular for every $t$. The crucial
fact, which helps us   avoid these difficulties  is that $\Bbb{H}^{3}_{1}$ has
dimension greater than  $G(T)$.

\begin{lemma} $\Bbb{H}^{3}_{t}$ is negatively curved for every $t$.

\end{lemma}

The proof of this statement is a standard application of the Gauss formula for
the curvature of a hypersurface.

\begin{theo} The map $U_{t}:  \widetilde{X} \longrightarrow \Bbb{H}^{3}_{t}$ is
pluriharmonic for every $t>0$.

\end{theo}

{\bf Proof.} Observe that the singular set of $G(T)$ consists of the centers of
the ideal hyperbolic polyhedra   and  the edges, which  come out of them.
Therefore we need to worry  only about ball  $b$ around a regular points $x_{0}
$ in $\widetilde{X}$ which maps to the ideal polyhedra on the boundary of
$\Bbb{H}^{3}_{t}$ which can be retracted to $G(T)$. We  consider 3 different
cases:

A) The map $U_{t}:  \widetilde{X} \longrightarrow \Bbb{H}^{3}_{t}$ maps  a ball
 $b$ around a regular point $x_{0}
$ in $\widetilde{X}$ to a face in the boundary of  $\Bbb{H}^{3}_{t}$ but far
from the edges . In this case the proof  of the theorem follows from Theorem  (
7.3) \cite{1}.

B) The  image of  a ball  $b$ around a regular point $x_{0}
$ in $\widetilde{X}$ under the map $U_{t}:  \widetilde{X} \longrightarrow {\Bbb
H}^{3}_{t}$   contains an edge comming out from the center of an ideal
hyperbolic polyhedra. In this situation
$\Bbb{H}^{3}_{t}$ looks localy like a product of a tree and $\Bbb{R}^{2}$.
Therefore Theorem  ( 7.2) \cite{1} applies  and we can find for every $t$ a
sequence of functions $ \psi_{i,t}$  which:

1)  vanish in a neighborhood  of the set $S_{t}$ of $U_{t}$;

2) tend to 1 on $X \backslash S_{t}$ and;

3) such that:

\[ \lim_{t\rightarrow \infty}{} \lim_{i\rightarrow 0}{\int_{X}^{}{\parallel
\bigtriangledown \bigtriangledown U_{t}}\parallel \psi_{i,t} d\mu}=0 .\]

 Following Theorem  ( 7.2) \cite{1} we get:

\begin{theo} Let $X$ be  a smooth projective variety,  $ \omega $ a parallel
p-form on $ \widetilde{X}$ and let  the  image of  a ball  $b$ around a regular
point $x_{0}
$ in $\widetilde{X}$ under the map $U_{t}:  \widetilde{X} \longrightarrow
\Bbb{H}^{3}_{t}$   contain an edge of an ideal hyperbolic polyhedra  but be far
from a vertex. Then in some small neighborhood around $x_{0}$  the form $\omega
\wedge dU$ satiesfies

\[ \delta_{\bigtriangledown}(\omega \wedge dU)\equiv 0.\]

\end{theo}

Here the notations are the same as in theorem 5.2.

Now in the same way as in Theorem  ( 7.3) \cite{1} we obtain:

\begin{theo} The map $U_{t}: \widetilde{X} \longrightarrow   \Bbb{H}^{3}_{t}$
is pluriharmonic.

\end{theo}

C) The  image of  a ball  $b$ around a regular point $x_{0}
$ in $\widetilde{X}$ under the map $U_{t}:  \widetilde{X} \longrightarrow
\Bbb{H}^{3}_{t}$   contains a translated center of an ideal hyperbolic
polyhedra. This case cannot be done as the previous two . $\Bbb{H}^{3}_{t}$
looks localy like a product  $P$ of an open polyhedra (namely we have taken the
botom face away)  and an open interval. We can approximate  $P$ by   products
$C_{t}$ of  intervals and cones. Following  the last example in the
Introduction of  \cite{1}, we can  also approximate the singular metric $g_{t}$
on $C_{t}$ by regular metrics $g_{ts}$ with $K\leq 0$.  This construction is
similar to the polyhedral immersion of a 2-disk into $\Bbb{H}^{4}$ descibed in
\cite{L}. To show that  the maps  $U_{ts} : \widetilde{X} \longrightarrow b$
are pluriharmonic for  every regular metric $g_{ts}$ on  $C_{t}$, we use the
same argument as in B). Namely we use Theorem  ( 7.2) \cite{1} to find for
every $t$ a sequence of functions $ \psi_{i,t}$.

 According to \cite{15} the limit of pluriharmonic maps is a pluriharmonic map
so we obtain a pluriharmonic map to $C_{t}$ with singular metric $g_{t}$ on it.
Now we approximate $P$ by the $C_{t}$, and using again  that limit of
pluriharmonic maps is a pluriharmonic map,  we get that the map  $U_{t}:
\widetilde{X} \longrightarrow  \Bbb{H}^{3}_{t}$ is pluriharmonic.

$\Box$

\begin{rem} The above argument works in the same way for a truncated  two
dimensional polyhedra  (\cite{32}). Part C) of the argument does not generalize
in higher dimensions.

\end{rem}

An easy corollary of the above theorem is the following:

\begin{corr} The map $U_{t}:  \widetilde{X} \longrightarrow G(T)$ is
pluriharmonic.

\end{corr}

  We are now ready to finish the proof of Theorem 6.1. After we get the
pluriharmonicity of the map $U:\widetilde{X}\longrightarrow G(T) $ we can
construct the spectral covering. As before the pluriharmonicity implies that if
we take the (0,1) part of the complexified differentials $du_{1}, du_{2}$ we
obtain holomorphic differentials over some   spectral covering.

The same argument as in the proof of Theorem 1.2  gives a complete
factorization  $ h:X \longrightarrow Y$. Let say that we are in case  A)   of
the Clemens-Lefschetz -Simpson theorem. Then we have  a local  isometry  $l:
\Bbb{R}^{2}\longrightarrow G(T)$  as in Lemma 4.2. We obtain $\Bbb{R}^{2}$ as
the real part of the map defined by the  integration of  the holomorphic forms
over the spectral covering $S$. The fact that  $ G(T)$ is strictly negatively
curved make  the existence of such a local  isometry impossible. Following
Lemma 4.2 we obtain in the same way a contradiction in the case when $dim_{\Bbb
{C}}  Y=2$.

According to the Clemens-Lefschetz -Simpson theorem there is a possibility of
an isometry  $l: \Bbb{R}^{1}\longrightarrow G(T)$. But  then if the image is
isometric to $ \Bbb{R}^{1}$ we will have subcomplexes   fixed under the action
of  $\pi_{1}(X) $ and this contradicts the surjectivity of
$\theta:\pi_{1}(X)\longrightarrow \pi_{1}(T)$.

Therefore  the  Clemens-Lefschetz -Simpson theorem  that there  is a
factorization through the fundamental group of an orbicurve.

$\Box$

A conjecture stated by J. Carlson and D. Toledo  \cite{31} says that if
$\Gamma$ is a K\"{a}hler group then $H^{2}(\Gamma)$ is nontrivial. As a simple
corollary of the above theorem 6.1 we have:

\begin{corr} If there exists a surjective homomorphism $\theta:\pi_{1}(X)
\longrightarrow \pi_{1}(T)$ then the above conjecture is true for $\pi_{1}(X)$.

\end{corr}

Inspired by \cite{22} we  describe a construction which we  use in the next
section. According  to the previous theorem, if we work with negatively curved
length spaces we get a nice description of the representations, namely they are
coming from the representations of the fundamental group of an orbicurve. It is
natural to try to hyperbolize the objects we are working with. This can be done
in the case of a building $B$. (Following \cite{23} we can replace  every
chamber in every apartment of the building by a negatively curved complex with
boundary the boundary of the chamber. )

According to  the result of R. Charney and M. Davis \cite{23} we obtain the
hyperbolized building  $HB$, which is strictly negatively curved. The
functoriality of the construction  of Charney and Davis  provides  us with an
action of $G$, the corresponding group of the initial building  $B$, over $HB$.
Unfortunately there is  one  disadvantage of the whole construction - $HB$ has
a nontrivial topology. So  we  can use the technique  developed above only on
representations $\varrho: \pi_{1}(X) \longrightarrow G$, for which we know
there exists a $G$- equivariant  continuous map
$U:\widetilde{X}{\longrightarrow}HB$.  Along the lines of the previous
corollary  we have:

\begin{corr} Every Zariski dense representation $\varrho: \pi_{1}(X)
\longrightarrow G$
for which there exists an $G$- equivariant  continuous  map
$U:\widetilde{X}{\longrightarrow}HB$, such that $G$ acts on $HB$ without a
fixed point, factors through the representation of the fundamental group of an
orbicurve.

\end {corr}

Let us  give some sufficient  conditions for the existence of  a $G$-
equivariant  continious map $U:\widetilde{X}{\longrightarrow}HB$.

\begin{theo}  There exists a  $G$- equivariant  continuous map
$U:\widetilde{X}{\longrightarrow}HB$ if  the group homomorphism
$\pi_{1}(X)\to Out (PI)$ induced by  $\varrho: \pi_{1}(X) \longrightarrow G$
can be lifted to a group
homomorphism $\pi_{1}(X)\to Aut(PI)$. Here $PI$ is the fundamental group of
$HB$ and $Out(PI)$  and $Aut(PI)$ denote the groups of
outer automorphisms and  automorphisms of $PI$ respectively.

\end{theo}

{\bf Proof.}

It is easy to see the  necessity of this condition  directly. It has
a natural interpretation in terms of gerbes (see \cite{34} and \cite{35}). Let
us explain the construction of the gerbes we  use in this situation. For this
we  need to choose a universal covering $\widetilde{HB}\to HB $ of
$HB=K(PI,1)$.

Then we have an exact sequence of groups:

\[ 1\longrightarrow PI \longrightarrow Homeo(\widetilde{HB}, HB)
\longrightarrow Homeo(HB) \longrightarrow 1 ,\]

where $Homeo(\widetilde{HB}, HB)$ is the group of homeomorphisms of
$\widetilde{HB}$ , which cover some homeomorphism of $HB$.

Since we have a principal bundle over $HB$  with group $Homeo(HB)$, there is
a gerbe whose objects are local liftings of the structure group to
$Homeo(\widetilde{HB})$, and whose arrows are isomorphisms of
$Homeo(\widetilde{HB})$ - bundles, which
induce the identity on the $Homeo(HB)$-bundle. The band of this gerbe is
locally isomorphic to the band associated to the constant sheaf of
$PI$ , but there is an outer twisting  represented by a class in
$H^{1}(X,Out(PI))$,  namely the class of the
homomorphism $\varrho: \pi_{1}(X) \longrightarrow G$. The obstruction to
realizing the given band as arising from some sheaf of groups over $X$  is
exactly  the obstruction
to lifting the above class to a class in $H^{1}(X , Aut(PI))$.

$\Box$

 Let us see how all this applies to the situation, where we let     $HB$  be
the hyperbolized $SL(3,\Bbb(Q)_{p})$ building. It easy to see that in this case
$PI=F_{SL(3,\Bbb(Q_{p}))}$, where $F_{SL(3,\Bbb(Q_{p}))}$  is  a free group
with $SL(3,\Bbb(Q_{p}))$ many  generators. $Out(F_{SL(3,\Bbb(Q_{p}))})$
contains $SL(3,\Bbb(Q_{p}))$. What we need  now, to claim the existence of a
$G$- equivariant  continuous map $U:\widetilde{X}{\longrightarrow}HB$, is the
splitting of the following exact sequence:

\[1\longrightarrow F_{SL(3,\Bbb(Q_{p}))} \longrightarrow
Aut(F_{SL(3,\Bbb(Q_{p}))} \longrightarrow SL(3,\Bbb(Q_{p})\longrightarrow 1,\]

over $T$, where $T$ is the set theoretic image of $SL(3,\Bbb(Q_{p}))$ in
$Out(F_{SL(3,\Bbb(Q_{p}))})$.

To have this splitting we need to make  $SL(3,\Bbb(Q_{p}))$ act on the graph we
have assigned  to $F_{SL(3,\Bbb(Q_{p}))}$ (the graph  in this case looks like
rose with one vertex and $SL(3,\Bbb(Q_{p})$ edges) while  fixing the vertex.
For this we need to make a  canonical choice of an initial point for $PI$  in
$HB$. A  sufficient  condition for this is, for example, if the image of
$\widetilde{X}$ under $U$ is  homotopy equivalent to a tree  which is  fixed
under the action of the image of $\varrho: \pi_{1}(X) \longrightarrow G$. If
this condition is satisfied,  we  change in a  canonical way  the generators of
$F_{SL(3,\Bbb(Q_{p}))}$  when we are change the initial point for
$F_{SL(3,\Bbb(Q_{p}))}$ in $HB$. One case when the above  condition  is
satisfied is when the whole image of $U$ is a tree  which is  fixed under the
action of the image of $\varrho: \pi_{1}(X) \longrightarrow G$.  This  is  an
equivariant  harmonic map $U$ having dimension 1. Unfortunately in general this
condition is hard to  check.

Using the technique of F. Paulin \cite{40} we prove:

\begin{theo} Let  X   be a smooth projective variety, and let  $\Gamma $ be  a
word hyperbolic group acting on the corresponding $\Bbb R$-tree without fixing
a vertex . Let  us assume that  $Out(\Gamma)$ is an infinite group and  let
$\rho:\pi_{1}(X) \longrightarrow \Gamma$ be a surjective homomorphism. Then
$\rho:\pi_{1}(X) \longrightarrow \Gamma$ factors through a representation of
the fundamental group of an orbicurve.
\end{theo}

{\bf Proof.} According to Paulin \cite{40}  $\Gamma$ has an fixed point free
action on $T$, an  $\Bbb R$ tree. Using a theorem of Reeps  \cite{40} we
conclude that this action is simplicial. We  obtain $T$ by deforming the Cayley
 graph of $\Gamma$ , using that $Out(\Gamma)$ is an infinite group.  From the
Schoen- Korevaar - Jost theorem we know that there exists a harmonic $
\pi_{1}(X)$-equivariant map to $T$. In the same way as in the proof of theorem
6.1 we show that $\rho: \pi_{1}(X) \longrightarrow \Gamma$ factors through a
representation of the fundamental group of an orbicurve. We  rule out the
possibility of part A) of the  Clemens-Lefschetz -Simpson theorem by using  the
fact that in the situation above we do not have parabolic action of $\Gamma$ on
infinity of $T$, namely  $\Gamma$ does not have fixed point of the   infinity
of $T$.

Another way to show that  $\rho(\pi_{1}(X)): \longrightarrow \Gamma$ factors
through a representation of the fundamental group of an orbicurve, is to
follow the proof of theorem 8.1 \cite{15}.

$\Box$

\begin{rem} We think that one can use the harmonic map technique to prove the
Morgan - Shalen conjecture for K\"{a}hler groups .The conjecture states
\cite{40}, that every group with free action on an $\Bbb R $ tree is a free
product of an abelian  and a surface group.

\end{rem}

\section{Some results about the integrality of representations}

 In this section we  discuss  a   partial verification of the following two
conjectures of Simpson \cite{11}. We are going to consider representations
$\varrho: \pi_{1}(X) \longrightarrow G$ , where $G$ is simple Lie group
defined over $k$, an  algebraically closed field of characteristic zero. In
this section we work only with $G=SL(n,k)$.

\begin{con} Let $\varrho: \pi_{1}(X) \longrightarrow G$ be  a rigid semisimple
representation. Then $\varrho \otimes k$   is a direct summand  over $k$ in the
monodromy representation of a motive (i.e. comes from geometry  ) over $X$.

\end{con}

\begin{con} Let $\varrho: \pi_{1}(X) \longrightarrow G$ be a rigid semisimple
representation. Then it is integral, in other words it is conjugate to a
representation whose matrix coefficients  are algebraic integers.

\end{con}

 (Note that the first Conjecture would imply the second.)

 Here rigid  representation means  that every representation which is nearby in
 the
  affine  variety of representations is conjugate to it. In the case of an
irreducible representation (the case we are working with), rigid is equivalent
to the fact that the corresponding point in the moduli space  of
representations is isolated.

The goal of this section is to verify the second  conjecture in some cases.

  We proceed following  Simpson \cite{3} . We assign to every Zariski dense
rigid representation  a new  Zariski dense rigid representation $\varrho:
\pi_{1}(X) \longrightarrow G^{1}$ over a local nonarchimedean field. The
procedure goes as follows. Observe that the  moduli space of representations is
defined over $\Bbb Q$, and since we are working with a rigid representation we
can find an isomorphic representation defined over $\bar{\Bbb Q}$.
  Let $E$ be a finite extension of  $\Bbb Q$ defined to  be the extension
which contains
 all coefficients of our representation, and let $O$ denote the ring of
integers in $E$.
  Let $E_{p}$ denote the field of fractions of the completion of $O$  in $p$ ,
for some prime $p$. Let $G^{1}$ be the new group
over  $E_{p}$ and use $\varrho$ again  for the representation $\varrho:
\pi_{1}(X) \longrightarrow G^{1}$. Since $E_{p}$ is a local field  then we are
in the situation of  the Theorem 4.2. Namely,
 using the Bruhat-Tits theory we can attach to $G^{1}$ a building $B$.
Following the Gromov and Schoen theory we can attach a   harmonic map to
$\varrho$. Therefore we get:

\begin{corr} Let $\varrho: \pi_{1}(X) \longrightarrow G$ be a Zariski dense
rigid representation. Then one of the following holds:

A) For every prime $p$ the image of $\varrho: \pi_{1}(X) \longrightarrow G^{1}$
is contained in some maximal compact subgroup in $G^{1}$.

B)For   some prime $p$ the image of $\varrho: \pi_{1}(X) \longrightarrow G^{1}$
is not  contained in some maximal compact subgroup in $G^{1}$. Then there
exist:

1)  a finite etale cover $X'$ of a blow up of $X$;

2) a smooth projective variety $Y$ of positive dimension    $l$ less or equal
to the rank $r$  of $G$ over $\Bbb{C}$;

3) and a  holomorphic map $h:X' \longrightarrow Y $

such that $\varrho:\pi_{1}(X') \longrightarrow G$ factors through a
representation of $\pi_{1}(Y) $.

\end{corr}

 In case A) we apply theorem os Simpson's \cite{C} to obtain that
$ \rho$ is an integer representation, namely a representation $\varrho:
\pi_{1}(X) \longrightarrow GL(r,{\Bbb {Z}})$, $r \geq n$.

Now we define special type of representations for which case B) of the above
theorem does not occur, namely conjecture 7.2 is satisfied. The idea is that
since curves do not have many rigid representations then if something factors
through a curve should not have either.

\begin{defi} Take  $\varrho: \pi_{1}(X) \longrightarrow G^{1}$ to be a Zariski
dense representation in $G^{1}$  which is assigned to $\varrho: \pi_{1}(X)
\longrightarrow G$. Define the corresponding   forms
 $de_{1}(z_{1}), \ldots , de_{1}(z_{1})$, and  the holomorphic one forms
$\alpha _{1}, \alpha _{2},..., \alpha _{n} \neq 0$ on $\widetilde{S}$ . If  all
of these forms span  a one dimensional subbundle  of the cotangent bundle of
$\widetilde{S}$  over a Zariski open set we  say that $\varrho: \pi_{1}(X)
\longrightarrow G$ is of dimension one.

  \end{defi}

\begin{theo} Let $\varrho: \pi_{1}(X) \longrightarrow G$ be  a rigid Zariski
dense representation of dimension  one. Then it is integral, in other words it
is conjugate to a representation whose matrix coefficients  are algebraic
integers.
   \end{theo}

 {\bf Proof.}

  The case of $G=SL(2,\Bbb{C})$ - representations was proven by Simpson
\cite{3} for any rigid  representation. Note that every Zariski dense
representation to  $SL(2,\Bbb{C})$  is of rank one.
For groups of higher ranks   we  use  the fact that our representation is of
dimension one. In this situation we can  use Theorem 1.1, namely, using the
Bruhat-Tits theory we can attach to $G^{1}$ a building $B$ and then attach,
following Gromov-Schoen theory, a harmonic map to $B$. Since our representation
is of dimension one, Theorem 1.2 gives us that there exists a finite
nonramified covering $X^{1} \longrightarrow X $ such that the fundamental group
of $X^{1}$ factors though the fundamental group of some smooth algebraic curve
$Y^{1}$. Now following \cite{3} we can actually factor the original
representation $\varrho: \pi_{1}(X) \longrightarrow G^{1}$ working with the
fundamental group of an orbicurve,   $\pi_{1}(Y,O)$, instead of with the
fundamental group of  $\pi_{1}(Y^{1})$.

 We proceed by applying a theorem of N. Katz.

\begin{theo} (Katz)   If $\varrho:\pi_{1}(Y,O)\longrightarrow G$ is
cohomologicaly rigid then it is  motivic.

\end{theo}

 To be able to apply the above theorem we need to show that rigidity in terms
of moduli spaces and the cohomological rigidity define the same objects. This
is the subject of the   next two lemmas.

Let us first give the precise definitions .

\begin{defi} Let $Y \cong  \Bbb{P}^{1}$. We  say that
$\varrho:\pi_{1}(Y,O)\longrightarrow G$ is a  physically rigid representation
if and only if  the conjugacy class  $\varrho:\pi_{1}(Y,O)\longrightarrow G$ is
  uniquely determined by its local monodronomies.

\end{defi}

 Let us introduce a new notion-the notion of cohomological rigidity.

\begin{defi} Let $Y \cong  \Bbb{P}^{1}$ and let
$\varrho:\pi_{1}(Y,O)\longrightarrow G$ be a   representation. Let  $F$  be the
bundle over $Y \cong  \Bbb{P}^{1}$  corresponding  to
$\varrho:\pi_{1}(Y,O)\longrightarrow G$.  We  say that
$\varrho:\pi_{1}(Y,O)\longrightarrow G$ is cohomologically rigid if and only if

\[\chi (X, j_{*}(End F))=2.\]

\end{defi}

For more detailed treatment of the notions of  physical  and cohomological
rigidity see \cite{36}.

Now we show that the notions of  rigidity in terms of moduli spaces and the
cohomological rigidity coincide.

\begin{lemma} Let $\varrho:\pi_{1}((Y,O)) \longrightarrow G$ be  a   Zariski
dense representation which is rigid in  the moduli sense. Then $Y \cong
\Bbb{P}^{1}$ with orbistructure in finitely many points.

\end{lemma}

{\bf Proof.} Let $a_{1},\ldots ,a_{g},b_{1},\ldots ,b_{g}, r_{1},\ldots ,r_{t}$
be  the generators of
$\pi_{1}((Y,O))$ .(Here $g$ is the genus of $Y$, and $t$ is the number of the
points with an orbistructure ). Consider the map:

\[ s:G_{1}\times\ldots \times G_{2g} \times C_{1} \times\ldots \times C_{t}
\longrightarrow G. \]
 where $G_{1},\ldots   , G_{2g}$ are $2g$ copies of $G$ and $C_{1},\ldots ,
C_{t}$ are  the conjugacy classes $SL(n,\Bbb{C})/Z(r_{1}),$ $\ldots ,$
$SL(n,\Bbb{C})/Z(r_{t})$.
Here the  $Z(r_{i})$ are the centralizers of the $\rho(r_{i})$ in
$SL(n,\Bbb{C})$. This map is defined to be just the multiplication of the
corresponding elements. The dimension of the maximal component of the  preimage
of 1 under this map gives an estimate for the dimension of the tangent space
to the moduli space of representations $\varrho:\pi_{1}(Y,O)\longrightarrow G$.
This dimension is equal to:

\[2g \dim G+\sum \dim SL(n,\Bbb{C})/Z(r_{i})- \dim G.\]

Since  $\varrho:\pi_{1}(Y,O)\longrightarrow G$ is a rigid representation, then
this dimension should be equal to  zero. But the number we have computed above
is always positive  ( $n\geq 3$ ) unless   $Y \cong  \Bbb{P}^{1}$.

$\Box$

The following  was shown in \cite{36}:

\begin{theo} (Katz) Let $Y \cong  \Bbb{P}^{1}$. Then the notions of physical
rigidity and cohomological rigidity over $Y$ are the same.

\end{theo}

The second lemma we need is the following:

\begin{lemma} Let $Y \cong  \Bbb{P}^{1}$. Then the notions of physical rigidity
and rigidity in terms of the moduli space of representations for an orbigroup
over  $Y$ are the same.

\end{lemma}

{\bf Proof.} Let us first show that physical rigidity implies moduli space
rigidity. Observe that $ r_{1},....,r_{t}$ are semisimple elements in $G$. We
need to show  that if  we have two representations of  $\pi_{1}(Y,O)$ that we
can deform one to   another, then these  two representations are conjugate,
providing  we know  they are physically rigid.  But since $ r_{1},....,r_{t}$
are semisimple,  if we deform them they need to stay in the same conjugacy
classes. Therefore they have the same local monodromies,  and  the physical
rigidity implies that they are actually conjugate.

Now we show that moduli space rigidity implies physical rigidity. Let us start
with a rigid, in terms of  moduli spaces, representation
$\varrho:\pi_{1}(Y,O)\longrightarrow G$.
We have already estimated  the dimension of the component of the moduli space
of representations which contains $\varrho$. Namely this dimension is less than
or equal to:

\[ \sum \dim SL(n,\Bbb{C})/Z(r_{i})- \dim G.\]

Observe that:

\[\sum_{i=1} ^{t} \dim Z(r_{i})=\left( \sum_{i=1} ^{t} \sum_{ j=1} ^{k_{i}}
n_{ij}^{2}\right)-1.\]

(Here $ k_{i}$ is the number of the different eigenvalues of $ \rho(r_{i})$ and
 $n_{ij}$ are the sizes of the corresponding  Jordan blocks.)

We also have:

\[\dim SL(n,\Bbb{C})/Z(r_{i})=(n^{2}-1)-(\left(  \sum
n_{ij}^{2}\right)-1)=n^{2}- \sum n_{ij}^{2}.\]

Since $\varrho:\pi_{1}(Y,O)\longrightarrow G$ is rigid in terms of the moduli
spaces we have that:

\[\dim SL(n,\Bbb{C})\geq  n^{2}-1=\sum_{i=1}^{t}\left\{  (n^{2}-1)-(\sum_{ j=1}
^{n} n_{ij}^{2})-1\right\}-(n^{2}-1).\]

This is equivalent to:

\begin{equation}
 2(n^{2}-1) \geq \sum _{i=1}^{t} (n^{2}-\sum_{ j=1} ^{n} n_{ij}^{2})
\end{equation}
But from the Grothendieck -Ogg-Shafarevich formula (see  \cite{36} ) we know
that:

\begin{equation}
\chi(X, j_{*}(End F)) = (2 - t) (n^{2}) +   \sum_{i=1, j=1} ^{t,n} n_{ij}^{2}.
\end{equation}
( Here $\chi(X, j_{*}(End F))$ is  the Euler characteristic of $ j_{*}(EndF)$
and $j$ is the embedding  $ j: \Bbb{P}^{1} -\; finite \; set \; of   \; points
\hookrightarrow  \Bbb{P}^{1}$.

Combining  (1) and (2)   we obtain:

\[\chi(X, j_{*}(End F)=\geq 2.\]

On the other hand  we have started with  an irreducible  representation
$\varrho:\pi_{1}(Y,O)\longrightarrow G$ we have:

\[\dim H^{0}(X, j_{*}(End F))= \dim H^{2}(X ,j_{*}(End F))=1\]

 Then $\dim H^{1}(X, j_{*}(End F))=0$ , which means that   the representation
$\varrho:\pi_{1}(Y,O)\longrightarrow G$ is cohomologicaly rigid. Now we are in
a position to apply theorem 7.2 and we conclude that
$\varrho:\pi_{1}(Y,O)\longrightarrow G$ is motivic.

But if $\varrho:\pi_{1}(Y,O)\longrightarrow G$ is motivic then the image of
$\varrho:\pi_{1}(Y,O)\longrightarrow G^{1}$ is in a bounded subgroup in
$G^{1}$. This can be seen as follows:

The fact $\varrho:\pi_{1}(Y,O)\longrightarrow G$ means that as representation
$\varrho$ extends from a discrete representation of the fundamental group of  a
complex variety extends to a representation of the profinite completion of the
fundamental group i.e. it extends to a representation of  $X$ as an algebraic
variety. But under continuous representation the compact  profinite completion
of the fundamental group goes to a compact image . So we get that the image of
$\varrho:\pi_{1}(Y,O)\longrightarrow G^{1}$ is  in a bounded subgroup in
$G^{1}$.

 Therefore the corresponding harmonic map is the constant map so the forms
$\alpha _{1}, \alpha _{2},..., \alpha _{n} \neq 0$ do not exists. This
contradicts our assumption that we have a factorization through an orbicurve.
So the second case does not occur.

$\Box$

In the process of the proof of this theorem we have obtained the following
corollary:

\begin{corr}  The preimage of 1 under the map defined above
\[s:C_{1}\times... \times C_{t}\longrightarrow G,\]
is connected.
\end{corr}

\begin{rem} The argument above depends a lot on the fact that we are working
with
$SL(n,\Bbb{C})$ (the dimensions of $SL(n,\Bbb{C})/Z(r_{i})$, a duality argument
which is hidden in the computation   of $\dim H^{2}(X ,j_{*}(End F))=1$), but
we believe that  it  is true for all  simple complex Lie  groups.

\end{rem}

Now we are going to extend  the field of application of the above theorem to
the quasiprojective case.

\begin{theo} Let $X$ be a quasiprojective variety and let $\varrho: \pi_{1}(X)
\longrightarrow G$ be  a rigid Zariski dense representation of dimension  one
and type B) with unipotent monodromy at infinity.

 Then $\varrho: \pi_{1}(X) \longrightarrow G$ is integral.

   \end{theo}

 {\bf Proof.} From  section 5 we know that we can define a spectral covering
for   quasiprojective variety. The proof of  theorem 7.4  just repeats   the
arguments in the proof of   theorem 7.2.

$\Box$

Now we are going to formulate 3   different sufficient conditions for which we
get a factorization through an orbicurve:

1)It follows from the work of Arapura, Bressler and Ramachandrachan  \cite{19}
that if  $H^{1}EL_{2}(\pi_{1}(X))$ \linebreak $\neq 0$ then every Zariski dense
representation $\varrho: \pi_{1}(X) \longrightarrow G$  factors through  a
representation of the fundamental group of an orbicurve.

2)It follows from \cite{20}   that if $H^{1}EL_{2}(S) \neq 0$, where $S$ is the
factorizing spectral covering, then every Zariski dense representation
$\varrho: \pi_{1}(X) \longrightarrow G$  factors through  a representation of
the fundamental group of an orbicurve.

3)It follows from \cite{2} that if  $Prym_{\sigma}(S,X)$ is nontrivial then
every Zariski dense representation $\varrho: \pi_{1}(X) \longrightarrow G$
factors through  a representation of the fundamental group of an orbicurve.

Therefore we arrive at:

\begin{corr} If $\varrho: \pi_{1}(X) \longrightarrow G$  is a rigid Zariski
dense representation and one of the above conditions is satisfied, then this
representation is integral.

\end{corr}

  There is one more case for which  we were able  to prove the fact that every
rigid representation is integral.

\begin{corr} Let  X be a smooth projective variety such that $\pi_{1}(X) $ is
a word hyperbolic group, acting  on the corresponding  $\Bbb R$-tree. Assume
that  $Out(\Gamma)$ is an infinite group. Let $\rho:\pi_{1}(X) \longrightarrow
SL(n,\Bbb C)$ be a rigid Zariski dense representation  such that
$Ker(\rho)\supseteq Ker(R)$. Then $\rho: \pi_{1}(X) \longrightarrow  SL(n,\Bbb
C)$ is integral.

Here $Ker(R)$ is the kernel of the map

\[R:\pi_{1}(X) \longrightarrow Isom(T),\]

where $Isom(T)$ is the group of isometries of the tree $T$,  defined in the
proof of theorem  6.7.

\end{corr}

{\bf Proof.} It follows  from theorems 6.1 ,  7.1  and 6.7.

$\Box$

In some cases  it is possible to prove  that every rigid representation is
motivic without assuming that the  harmonic map to the building  $B$ is of
dimension one. It follows from Corrolary 6.2 that:

\begin {corr} Every rigid Zariski dense representation $\varrho: \pi_{1}(X)
\longrightarrow G$ for which there exists an $G$- equivariant nonconstant
continuous  map $U:\widetilde{X}{\longrightarrow}HB$, is integral.

\end{corr}

{\bf Proof.} If the action of $G^{1}$ fixes a point in  $HB$ then  the image of
$\varrho: \pi_{1}(X) \longrightarrow G$ is contained in a maximal compact
subgroup of $G$. Using again the theorem of Baas  we get that our
representation is integral. Now if  we have an action of $G$ without a fixed
point we are in a position to apply  Corollary 6.2.
In this case it is very easy to prove (see \cite{16}) the pluriharmonicity of
$U$  since we keep almost the same metric on the boundary of every chamber.
Then  we get a  spectral covering as in Theorem 7.1. Also  the spectral
covering, which we  obtain has an  involution on it. Following \cite{3} and
\cite{4} we obtain an involution which acts on the curve $Y$ as well. We  mod
out the spectral covering $S$ and the curve  $Y$ by $\Bbb{Z}_{2}$ and we  have
a holomorphic map between the quotients. The rest of the proof is the same as
in the main theorem. We finish it   by applying theorem 7.3.

$\Box$

The above corollary actually implies that if there exists an $\rho$-
equivariant nonconstant  continuous  map $U:\widetilde{X}{\longrightarrow}HB$
and $ \rho$ is a rigid representation then the action of $\rho$ on $B$ should
have a fixed point.

We think of the above corollary as a generalization of the tree case, since the
tree is negatively curved. Unfortunately as we saw in Theorem 6.3 it is in
general  impossible to  make $G$  act on simply connected negatively curved
spaces. So one should try another   approach for the  general Simpson
conjecture. It seems  one way of doing that would be to  generalize the result
of Katz to higher dimensions. While at the moment the Conjecture 7.1 and
Conjecture 7.2 look out of reach in their complete generality the following
conjecture looks more doable at this moment.

\begin{con} Every rigid Zariski dense representation $\varrho: \pi_{1}(X)
\longrightarrow G$ of a group $ \pi_{1}(X)$ with the property $T$ is
integrable.
\end{con}

\section{The Shafarevich conjecture}

Inspired by the results in \cite{30} we make in this section a connection
between the
Shafarevich map and the results of this paper.  We recall the definition of
the Shafarevich map \cite{30} , and describe some of its properties.

\begin{defi} Let $X$ be a smooth projective variety. Then we call a  rational
map

\[Sh:X ----\rightarrow Sh(X),\]
the Shafarevich map of $X$ iff:

1) $Sh(X)$ is a normal projective variety;

2)The rational map  $Sh:X \longrightarrow Sh(X)$ has connected fibers and;

3) there are countable many closed subvarieties $D_{i}$ in $X$ ($D_{i}\neq X$)
such that for every irreducible subvariety $Z$ in $X$,  not  in the union of
the $D_{i}$
we have:

\[Sh(Z)= {\rm point \; iff} im[\pi_{1}(Z')\longrightarrow \pi_{1}(X)]  {\rm is
\; finite}.\]
\end{defi}

It is easy to see that if the rational map  $Sh:X \longrightarrow Sh(X)$
exists, it is unique up to birational equivalence.

One  can give a relative version of this definition with respect to a normal
subgroup $H$ of  $\pi_{1}(X)$ .To do that one needs to   require in part 3) of
the definition above that  $im[\pi_{1}(Z')\longrightarrow \pi_{1}(X)] \cap H $
 has a finite index in  $im[\pi_{1}(Z')\longrightarrow \pi_{1}(X)]$. In this
case we  write $Sh^{H}:X \longrightarrow Sh^{H}(X)$ for the corresponding
rational map.

The following theorem is proven in \cite{30}.

\begin{theo}( J\'anos Koll\'ar) The rational maps $Sh:X \longrightarrow Sh(X)$
and $Sh^{H}:X ---\longrightarrow Sh^{H}(X)$ exist.

\end{theo}

The properties of the Shafarevich map we have described above provide us with
enough information to prove the following theorem:

\begin{theo} Let X be a smooth projective variety  which has type B (from
Theorem 1.1) Zariski dense representation of its fundamental group to a Lie
group $G$ defined over a field with  discrete valuation. Then:

\[dim_{\Bbb{C}}Sh^{H}(X)=dim_{\Bbb{C}}(Y) \leq rank_{\Bbb{C}}(G),\]

where Y is the variety defined in Theorem 1.1 and $H=Ker (\varrho) $. Moreover
we have that some  finite nonramified coverings of $Y$ and $ Sh(X(\Gamma)$
(defined below) are birationally isomorphic.

\end{theo}

{\bf Proof.} From the construction of the  holomorphic map $X  \longrightarrow
Y$ in the proof of Theorem 1.2  we know that if $Z$ is the general fiber of
this map we have that the intersection of the subgroups
$R=im[\pi_{1}(Z)\longrightarrow \pi_{1}(X)]$  and   $H=Ker (\varrho) $ is a
subgroup of  finite index in $R$. This implies that $Sh^{H}:X \longrightarrow
Sh^{H}(X)$ factors through $X \longrightarrow Y$ and we get:

\[ rank_{\Bbb{C}}(G) \geq dim_{\Bbb{C}} Y \geq dim _{\Bbb{C}} Sh^{H}(X).\]

To get the inverse inequality we argue as follows:

Let us first do two reductions.

Denote by $L$ the intersection of all subgoups of finite index in $\pi_{1}(X)$.
 According to remark 4.4 the conclusions of the main theorem are still true if
we work with $\pi_{1}(X) \not L$. Therefore we can assume that

\[ H=\bar{H}. \]

Define now $R=im[\pi_{1}(Z)\longrightarrow \pi_{1}(X)]$.

The second reduction is to show that we can always assume that  $R$ is
contained in $H$. If  $R$ is not  contained in $H$ then
$R=im[\pi_{1}(Z)\longrightarrow \pi_{1}(X)]$    intersects  $H=Ker (\varrho) $
in a subgroup of   finite index in $im[\pi_{1}(Z)\longrightarrow \pi_{1}(X)]$.
Then let the image of this finite subgroup  be the finite subgroup $ \Pi$ in $
\pi_{1}(X)/H $. Now we do the following  modification:

Observe   that  $ \pi_{1}(X)/H $ is a residually finite group due to the fact
that $G$ is an affine algebraic group. Therefore since we have mod out $
\pi_{1}(X)$ by $L$ and $H=\bar{H} $ we can find a   subgroup of finite index
$\Theta $ in  $ H=\bar{H}/H $ which does not intersect $ \Pi$. We conclude that
for the covering , which corresponds to the group  $\Theta $, $R$ is containe
in $H$.

So from now on we assume that $R$ is contained in $H$.

Applying  4.8 \cite{30} we see that there exists some subgroup $ \Gamma$  of
finite index in $ \pi_{1}(X)$ such that for the finite nonramified covering
$X(\Gamma)$ we have a complete factorization through $Sh^{H}(X(\Gamma))$.

The construction goes as follows. We have already  constructed $X(\Gamma)$ in
section 4.  According to 4.8 there exists a   map:

\[f_{*}:  \pi_{1}(X(\Gamma)) \longrightarrow \pi_{1}(Sh^{H}(X(\Gamma))) \]

such that

\[Ker f_{*} \subseteq  \bar{R} \subseteq \bar{H} =H.\]

But $Y$ was defined as the variety of minimal dimension, through which we have
a factorization. Therefore:

\[ dim_{\Bbb{C}}Sh^{H}(X) \geq dim_{\Bbb{C}} Y.\]

Combining this with the previous inequality  we get:

\[dim_{\Bbb{C}}Sh^{H}(X)=dim_{\Bbb{C}}(Y) \leq rank_{\Bbb{C}}(G).\]

$\Box$

We prove now a weak version of the Shafarevich conjecture.
The Shafarevich conjecture says that for every smooth projective variety $X$
there exists a there exists a Stein manifold ${\bf Sh} (\widetilde{X} )$ and a
proper map with connected fibers $ {\bf Sh}: \widetilde{X}\longrightarrow {\bf
Sh}(\widetilde{X})$.

An easy consequence of previous theorem is the following theorem:

\begin{theo}Let X be a smooth projective variety  which has a  type B (from
Theorem 1.2) Zariski dense representation $\rho$ of its fundamental group  to a
Lie  group $G$ defined over a  local field $K$. Define $H=Ker(\rho)$. Then $Y$
is isomorphic to $Sh^{H}(X')$ and the map

\[{\bf Sh}^{H}:X'\longrightarrow {\bf Sh}^{H}(X'),\]

is a morphism. Here $ X'$ and $Y$ are the same as in Theorem 1.2.

\end{theo}

{\bf Proof.}Without lost of generality we can work with the case when
$H=Ker(\rho)$ is finite.

It is clear that if $Z$ is a subvariety in $X$ and the fundamental group of $Z$
is finite then
$Z$ goes to a point in $B$, where $B$ is the corresponding building. This
follows from a very general theorem of M. Bridson saying that every finite
group, which acts on a CAT space has a fixed point. Therefore the corresponding
harmonic maps are just constants. We need to show that if $Z$ is a subvariety
of $X$ with a finite fundamental group then then it goes to a point in $Y$.

According to Gromov and Schoen the map $U:\widetilde{X}\longrightarrow B$ is
essentially regular, namely the intrinsic derivative exists everywhere (see
\cite{1}). We claim that $U^{*}\alpha_{i}$ is zero on $Z$. On the intersection
of the  smooth locus  of $U$ and $Z$  this is obviously true. For the singular
points we use the same argument as in the proof of    Lemma 3.1, namely the
existence of the intrinsic derivative. This existence  implies existence of a
kernel of $U^{*}$ in   $T^{*}X'$ and we conclude  $U^{*}\alpha_{i}$  is zero on
$ Z$. This implies that $Z$ will go to a point in $Y$  even if $Z$ it is
contained in  the singular set of $U$.

Another way to show that $Z$, a subvariety of $X$ with a finite fundamental
group  goes to a point in $Y$ is the following. Using theorem 6.3 \cite{1} we
obtain the following estimate, which  follows from the existence of the
intrinsic derivative

\[\lim_{x\rightarrow S}{\|dU(x)\|}=0,\]

where $S$ is the singular set of $U$.

According  to \cite{1} $U(S)$ is contained in the closed faces  of the
simplices of the highest dimension. Now using the above estimate and
approximating
$U(S)$ in the normal directions of the faces containing it
we get that $U^{*}\alpha_{i}$  is zero on $ Z$.
\hfill $\Box$

 Therefore  we obtain the following :

\begin{corr} If the conditions of the above theorem  we have $H=ker \rho$ is a
finite group then the Shafarevich-Koll\'ar conjecture   is true.
\end{corr}

{\bf Proof.} According to the previous theorem the map

\[{\bf Sh}:X'\longrightarrow {\bf Sh}(X'),\]

is a morphism. Consider  the morphism

\[h:X \longrightarrow Y \]

 from the proof of the main theorem in section 4.

 Assume that there is a subvariety $Z$ in $X$ with the property that
$im[\pi_{1}(Z)\longrightarrow \pi_{1}(X)]$ is finite and it does not go to a
point in $Y$. But then pulling it back to $X' $ we get a subvariety $Z^{'}$ in
$ X'$ such  that $ \pi_{1}(Z^{'})$ is a finite group due to a theorem of
Campana \cite{31}. Therefore
$im[\pi_{1}(Z^{'})\longrightarrow \pi_{1}(X')]$ is a finite group and $Z^{'}$
goes to a point in $ {\bf Sh}(X')$. Now using the map  $ {\bf Sh}(X')
\longrightarrow  Y $ we obtain a contradiction since the image of  $Z^{'}$ in
$Y$ is the same as the image of  $Z$ in $Y$. Therefore  $Y={\bf Sh}(X)$.
\hfill $\Box$

The  Shafarevich-Koll\'ar conjecture in general is connected with another
question of Gromov:

Can we find a faithful  discrete cocompact action of every K\"{a}hler word
hyperbolic group  to a  space with $K \leq 0$ ?

We believe that  the Shafarevich-Koll\'ar conjecture   follows for every
K\"{a}hler group for which this could be done.

Now we show two examples for which the above condition is satisfied and we can
actually prove the Shafarevich conjecture.

\begin{corr} Let $X$ is a compact K\"{a}hler  manifold and $\pi_{1}(X)$ is an
amalgamated product of two groups and $\pi_{1}(X)$ is word hyperbolic. Then the
 Shafarevich-Koll\'ar conjecture is true for $X$.

\end{corr}

{\bf Proof.} According to \cite{1} we have a  faithful discrete cocompact
action of $\pi_{1}(X)$ on a tree. The rest is just repeating the argument from
the proof of  Theorem 8.3.
\hfill $\Box$

In the same way   Theorem 6.1  implies :

\begin{corr} Let  $X$ be  a compact K\"{a}hler  manifold and $\pi_{1}(X)$ has
an a faithful discrete cocompact action on a two dimensional complexes  defined
in section 6.   Then the  Shafarevich-Koll\'ar conjecture is true for $X$.
\end{corr}

Now we  discuss some possible applications of the Shafarevich-Koll\'ar
conjecture.

 One possible application of the above statement is to try to answer the
following   question of M. Ramachandran \cite{M}.

\begin{con}( M. Ramachandran) Let X be a smooth projective variety  which has a
 type B (from Theorem 1.2) Zariski dense representation $\rho$ of its
fundamental group  to some  Lie  group $G$ defined over a  local field $K$.
Define $H=Ker(\rho)$ and assume $H$ is a finite group.
Then one of the following holds:

1) The universal cover of $X$ satisfies the Bochner - Hartogs property, namely
for every $\alpha \in A^{0,1}_{c}(\widetilde{X})$ satisfying

\[\bar{\partial} \alpha =0 ,\]

there exists $u \in  C^{\infty}_{c}(\widetilde{X})$   such that

\[ \bar{\partial} u = \alpha.\]

2) $\pi_{1}(X)$ is comesurable to the fundamental group of a compact Riemann
surface.

\end{con}

Here  $  A^{0,1}_{c}(\widetilde{X})$ is the space of compactly supported smooth
(0,1) forms on $\widetilde{X}$  and $ C^{\infty}_{c}(\widetilde{X})$ is the
space of  all compactly supported smooth functions  on $\widetilde{X}$ with
complex values.

\begin{rem} M. Ramachandran actually stated the conjecture in much bigger
generality namely for every  smooth projective variety.

\end{rem}

\bigskip

\noindent

We hope to answer this question by using the maps

\[U:\widetilde{X} \longrightarrow B\]

and

\[{\bf Sh}^{H}:X'\longrightarrow {\bf Sh}^{H}(X').\]

to produce enough plurisubharmonic  functions on $\widetilde{X}$.

We study now some questions connected with the topological nature of  the map

\[{\bf Sh}^{H}:X'\longrightarrow {\bf Sh}^{H}(X').\]

We give a  partial answer to some questions posed by J\'anos Koll\'ar
 at the end of  \cite{27} . These questions are about how  much the Shafarevich
variety  ${\bf Sh}(X)$ depends on the original variety $X$ e.g.  does a
deformation of $X$  give a deformation of ${\bf Sh}(X)$ an so on.

\begin{corr} Let  $X$  be a smooth projective  variety which has nonrigid
finite kernel Zariski dense representation in $SL(2,\Bbb{C})$.  Let
$X_{t}$ , $t \in $ to the unit disc $\Delta$ , be a holomorphic deformation of
$X$. Then ${\bf Sh}(X_{t})$ is a holomorphic deformation of ${\bf Sh}(X)$.
\end{corr}
{\bf Proof.} The moduli space of $\lambda$ - connections (see \cite{35})
deforms  together  with $X$ . The  moduli space of $\lambda$-connections
itself is a deformation of the moduli space of  $SL(2,\Bbb{C})$ representations
to the moduli space of $SL(2,\Bbb{C})$ Higgs bundles. Since we have fixed the
representation we get a  constant holomorphic section over  $\Delta\times
\Bbb{P}_{1}$, which specializes at $0$ of $\Bbb{P}_{1}$ to an Higgs bundle. As
we said the moduli space of $SL(2,\Bbb{C})$ Higgs bundles is a deformation of
the moduli space $\lambda$ - connections  so we get a  holomorphic section  to
the family of  moduli spaces of $SL(2,\Bbb{C})$ Higgs bundles. over $ \Delta$.
Also we have a holomorphic map from the   moduli space of $SL(2,\Bbb{C})$ Higgs
bundles to the  $Hilb (T^{*}(X))$, which corresponds to every Higgs bundle a
spectral covering (see \cite{3}).  We get this way a holomorphic  deformation
of the corresponding spectral coverings and as a consequence   a family  of
curves $C_{t}$ through which according to our main theorem  the representation
factors.

But since this is an finite kernel Zariski dense representation we
have ${\bf Sh}(X_{t}(\Gamma))=C_{t}(\Gamma)$, where  $C_{t}(\Gamma)$ is the
spectral covering of $X_{t}(\Gamma)$. Observe that the map ${\bf
Sh}(X_{t}(\Gamma))\longrightarrow {\bf Sh}(X_{t})$ has the same Galois group as
 $X_{t}(\Gamma) \longrightarrow X_{t}$.  Therefore we  can mod out by this
group and  get that ${\bf Sh}(X_{t})$ deforms itself. \hfill $\Box$

Even more is true:

\begin{corr} Let $X$ be a   smooth projective varietiy which has nonrigid
finite kernel  Zariski dense representation  $\rho:\pi_{1}(X) \longrightarrow
SL(2,\Bbb{C})$. Let $Y$ be  another smooth projective varieties which is
homeomorphic to  $X$. Then ${\bf Sh}(X)$ is  homeomorphic to $Sh(Y)$.

\end{corr}
{\bf Proof.} Due to the fact that we have a finite kernel Zariski dense
representation we conclude that the $\pi_{1}(X)$ is infinite and therefore
${\bf Sh}(X(\Gamma))$ is a curve. Therefore   $\pi_{1}(X)= \pi_{1}(Y)$ and let
say that $C_{X(\Gamma)}={\bf Sh}(X(\Gamma))$ and $C_{Y(\Gamma)}={\bf
Sh}(Y(\Gamma))$
are the curves through which our representation factors. The only thind we need
to show is that  $g(C_{X})=g( C_{Y})$, where $g(C_{X})$ is the genus of
$C_{X}$. According to a theorem of Siu ( see for example \cite{5} ) if we have
a map from $\pi_{1}(X)$ to $\pi_{1}(C_{Y})$  we have a holomorphic map
$h:X(\Gamma)\longrightarrow C^{'}_{Y}$ to some other curve $ C^{'}_{Y}$. Due to
the fact that $C_{X}={\bf Sh}(X(\Gamma))$ we get that $g(
C_{Y})=g(C^{'}_{Y})\geq g(C_{X})$. Now to finish the proof we change the places
of $X$ and $Y$ in  the above argument.

\hfill $\Box$

\begin{rem} The above corollary  works for $X_{t}$ homotopy equivalent to  $X$
as well as for ${\bf Sh}_{K}(X)$ instead of  ${\bf Sh}(X)$
, where $K$ is the kernel of $ \rho$ and $K$ is an infinite group.

\end{rem}

The following two corollaries are straightforward  consequences from section 7.

\begin{corr} The corollary above is  true  if instead of
a Zariski dense representation in $SL(2,\Bbb{C})$ we require the existence of a
finite kernel surjective  homomorphism of the  fundamental group of $X$  to the
fundamental group of a 2 dimensional negatively curved complex of groups (see
paragraph 6).
\end{corr}

\begin{corr} The corollary above is  true if instead of
a Zariski dense representation in $SL(2,\Bbb{C})$ we require the existence of a
nonrigid Zariski dense representation  of dimension  1 (see paragraph 7) to any
simple complex  Lie group.

\end{corr}

In the same way we treat the case when $Y$ is a surface.

\begin{theo} Let  $X$  be a smooth projective  variety which has  a finite
kernel Zariski dense representation $\rho:\pi_{1}(X) \longrightarrow G $ of
type B  in $G$, where $G$ is a simple Lie group such that $rank_{\Bbb{C}} G
=2$. Let
$X_{t}$ be a holomorphic deformation of $X$. Then ${\bf Sh}(X_{t}(\Gamma))$ is
a holomorphic deformation of ${\bf Sh}(X(\Gamma))$.

\end{theo}
{\bf Proof.} The assumptions of the theorem imply that we have a factorization
through the representation of the fundamental group of $Y$, where $Y$ is either
an algebraic curve or an  algebraic surface.  We have considered the first case
in the previous corollary.

Let us consider now the case $dim{\bf Sh}(X(\Gamma))=2 $.
Observe that ${\bf Sh}(X_{t}(\Gamma))$ is defined only up to birational
isomorphism. Therefore  when we say that  ${\bf Sh}(X_{t}(\Gamma))$ is a
holomorphic deformation of ${\bf Sh}(X(\Gamma))$,  we mean that the minimal
model of  ${\bf Sh}(X_{t}(\Gamma))$ is a holomorphic deformation of the minimal
model of  ${\bf Sh}(X(\Gamma))$. Since  $X_{t}$ is a holomorphic deformation of
$X$   and we work with the same representation    in the same way as in the
previous corollary, we can show that  the corresponding  spectral coverings of
$X_{t}$  and $X$ are holomorphic deformations of each other. What we need to
show is that if $Y$ is a  holomorphic deformation of $Y_{t}$ then the minimal
model of  $Y$ is a holomorphic deformation of $Y_{t}$. But in case of surfaces
with $ \kappa(Y) \geq 0 $ (this is exactly our case since $ {\bf
Sh}(X(\Gamma))$ has a large fundamental group) this follows from a theorem of
Iitaka. The fact that $Y$ and ${\bf Sh}(X(\Gamma))$ are  isomorphic proves the
theorem.\hfill $\Box$

 The following  nonvanishing theorem is a consequence of (3).

\begin{theo}

Let $X $ be a smooth fourfold of general type with a  type B  (from Theorem
1.1)  Zariski dense representation of  its fundamental group $\rho: \pi_{1}(X)
\longrightarrow G $. Here $G$ is a  rank 2 or 3  simple Lie group over field
with a discrete valuation. Then either:

 1) $\rho: \pi_{1}(X) \longrightarrow G $ factors through a representation of
the fundamental group of an orbicurve or;

2)  $P_{n}$=$H^{0} (X,  nK_{X}) $  for $n\geq4$ is not zero.

\end{theo}

{\bf Proof.} According to  Theorem 1.2 our representation factors through a
representation of the fundamental group of some variety $Y$ of dimension less
than or equal to 3. If this $dimY= 1 $ then we are in part 1) of the theorem.

Consider the case where the dimension of $Y$ is at least 2. The previous
theorem tells us that  $Y$ is birational to $Sh^{H}(X)$. Following  J\'{a}nos
Koll\'{a}r \cite{30}
(4.5 and 5.8) we can work with a morphism $X^{1} \longrightarrow S $, where
$X^{1}$ is a nonramified finite covering of $X$ and $S$,  a smooth variety
with generically large fundamental group.
It  follows from \cite{30}  that in this case nonvanishing theorems for $X^{1}$
imply nonvanishing theorems for $X$. We need to consider the following cases:

1) $dimS=3$   and $S$ is of general type. By  \cite{30} (10.1 ) we have that
$P_{m}(S)\geq1$ for $m\geq2$. Observe also that the fibers $X_{s}$ of the map
$X^{1}\longrightarrow S$ are curves of general type and therefore we have
$P_{1}(X_{s})\geq 1$.   Using \cite{30} (10.4)  we obtain:

\[ P_{n}(X) \geq P_{n-2}(S) \geq 1,  \; {\rm for}   \; n\geq 4.\]

2) $dimS=3$   and $S$ is abelian variety.  Then we can  apply \cite{30} (8.10
), which  gives a strong nonvanishing. Strong nonvanishing is equivalent to:

\[ h^{0}(X^{1},K_{X^{1}} \otimes D)\neq 0 \;  {\rm for} \;  D \;  {\rm big  \;
divisor \; on} \; X^{1} \; {\rm and} \; X^{1} {\rm big  \; birational \; to }
X. \]
In our situation we can make $3K_{X} = D $ since $X$ is of general type and
therefore:

\[P_{n}(X) \geq 1   \; {\rm  for}  \; n\geq 4. \]

  This  argument actually applies in the previous case too. We need  only that
the fundamental group of $S$ is  generically large.

3) $dimS=2$   and $S$ is abelian. We know that the fiber $X_{s}$ of
$X^{1}\longrightarrow S$ is a surface of general type. If we have a complete
factorization $X^{1}\longrightarrow S$ of the representation $\varrho$ we get a
contradiction, since $\varrho$ is a Zariski dense representation and we cannot
have a Zariski dense representation of a free abelian group in a simple Lie
group $G$. Therefore we need to work with a finite ramified cover of  $S$ - a
surface  of general type with generically large fundamental group, which we
will also denote by the  letter $S$. We are in a position to apply  \cite{30}
(10.4).

We have:

\[P_{n}(X_{s}) \geq 1 \; {\rm for} \; n\geq 2,\]

and also

\[ P_{n}(S) \geq 1 \; {\rm for}  \; n\geq 2.\]

Therefore we obtain from  \cite{30} (10.4):

\[ P_{n}(X) \geq P_{n-2}(S)\geq 1 \;   {\rm for} \; n\geq 4.\]

4) $dimS=2$   and $S$ is of general type. Then the   fiber  $X_{s}$ of
$X^{1}\longrightarrow S$ is also a  surface of general type. In this case we
argue as before,  again  using \cite{30} (10.4).

\begin{rem} Obviously  the parts  1)  and  2)  are not mutually exclusive.

\end{rem}

$\Box$

In the last application we are going to be concerned with the nonexistence of
some representations  for a special class of algebraic varieties. If we think
of $Hom (\pi_{1}(X), G)$ as the first nonabelian cohomology group, then the
nonexistence of certain representations is a sort of vanishing theorems.

\begin{defi}  We  say that an  algebraic variety   $X$ has a generically  large
fundamental group iff  for every subvariety  $Z$ of $X$ outside some union of
divisors $ \bigcup{D_{i}}$ we have that:

\[im[\pi_{1}(Z')\longrightarrow \pi_{1}(X)]\]

is an infinite group.

\end{defi}

\begin{corr} Let $X$ be  a variety  with a generically  large fundamental
group. Then:

1)A finite nonramified covering of  $X$ is birationaly isomorphic to $Sh(X)$
and

2) $X$  does not have Zariski dense,  finite kernel representations in a
complex Lie group $G$ if $dim_{\Bbb{C}}(X)>rank _{\Bbb{C}}(G)$.

\end{corr}

{\bf Proof.} Part 1) of this corollary follows immediately from the definition
of the Shafarevich variety.  Part 2) of this corollary follows  from theorem
8.3 since:

\[dim_{\Bbb{C}}X=dim_{\Bbb{C}}Sh(X)=dim_{\Bbb{C}}(Y) \leq rank_{\Bbb{C}}(G),\]

and this contradicts our assumption that:

\[dim_{\Bbb{C}}(X)>rank _{\Bbb{C}}(G).\]

$\Box$

\section{Final remarks}

By definition all of our spectral coverings are zero shemes of  sections in the
cotangent bundle of $X$. Then obviously they can be deformed to the zero
section. The zero  section taken with some multiplicity (in the  case of
buildings  this multiplicity was equal to  $w$) corresponds to a harmonic map
to a point in the Euclidean building $B$ corresponding to $G$. Therefore  we
have an action of  $\pi_{1}(X)$ on $B$ which fixes  a point. In all of our
arguments before we have excluded the case where $\pi_{1}(X)$ fixes a point  of
$B$ b. But if the fixed  point $P$  is not a point at infinity of $B$  we
assign to the action of $\pi_{1}(X)$ the constant equivariant harmonic map
$\widetilde{X }\longrightarrow P$. Obviously all holomorphic differentials in
this case are  equal to zero and the  spectral covering $S$  is just $X$, the
zero  section of the cotangent bundle of $X$. Therefore every nonrigid subgroup
 $G(B)$  (the group of  isometries of  $B$)  which cannot be deformed to
another subgroup of  $G(B)$ fixing a point of $B$, cannot be the fundamental
group of a smooth projective variety. Of course one needs to find  (if
possible) the proper notion of deformation. It is clear that such a statement
might then be true  in  much more general situations, namely for nonrigid
subgroups of the group of isometries of  quite general nonpositively  curved
length spaces.

Another direction one can try to apply the techniques developed in this paper
is to study the p-adic uniformization defined in \cite{13}, \cite{14},
\cite{15}.

Using the techniques described above one can try to study the discrete
subgroups of the infinite dimensional Lie groups $SDiff (D)$, where $D$ is a
Riemannian domain.
It follows from \cite{27}, \cite{28} and  \cite{29} that these groups are
nonpositively curved in some sense if  $D=T^{2} $, $D=T^{n} $ or  $D=S^{2} $.
Then using the harmonic map technique one obtains  some finiteness results.
These remarks  are going to  be  an object of future considerations.

\end{document}